\documentclass[aop,citesort,MSNbibl,dvips]{arximspdf}
\usepackage{pstricks}
\usepackage{pst-tree}
\usepackage{pst-node}
\usepackage{graphicx}

\doi{10.1214/09-AOP500}
\volume{38}
\issue{2}
\pubyear{2010}
\firstpage{532}
\lastpage{569}

\makeatletter

\makeatletter
\let\color@gray\@gobbletwo
\makeatother

\newproclaim{assumption}{Assumption}
\newtheorem{theorem}{Theorem}
\newproclaim{remark}{Remark}
\newtheorem{lemma}{Lemma}

\makeatother

\begin{document}
\begin{frontmatter}

\title{Taylor expansions of solutions of stochastic partial
differential equations with additive noise\protect\thanksref{T1}}
\runtitle{Taylor expansions for SPDE's}

\thankstext{T1}{Supported by the DFG project ``Pathwise numerics and dynamics
of stochastic evolution equations.''}

\begin{aug}
\author[A]{\fnms{Arnulf} \snm{Jentzen}\corref{}\ead[label=e1]{jentzen@math.uni-frankfurt.de}} and
\author[A]{\fnms{Peter} \snm{Kloeden}\ead[label=e2]{kloeden@math.uni-frankfurt.de}}
\runauthor{A. Jentzen and P. Kloeden}
\affiliation{Johann Wolfgang Goethe University, Frankfurt am Main}
\address[A]{Institute of Mathematics\\
Johann Wolfgang Goethe University\\
Robert-Mayer-Str. 6-10\\
D-60325 Frankfurt am Main\\
Germany\\
\printead{e1}\\
\phantom{E-mail: }\printead*{e2}} %adresu isvedimo komanda gale!
\end{aug}

% HISTORY:
\received{\smonth{9} \syear{2008}}
\revised{\smonth{7} \syear{2009}}

% ABSTRACT
%
\begin{abstract}
The solution of a parabolic stochastic partial differential equation
(SPDE) driven by an infinite-dimensional Brownian motion is in general
not a semi-martingale anymore and does in general not satisfy an
It\^{o} formula like the solution of a finite-dimensional stochastic
ordinary differential equation (SODE). In particular, it is not
possible to derive stochastic Taylor expansions as for the solution of
a SODE using an iterated application of the It\^{o} formula.
Consequently, until recently, only low order numerical approximation
results for such a SPDE have been available. Here, the fact that the
solution of a SPDE driven by additive noise can be interpreted in the
mild sense with integrals involving the exponential of the dominant
linear operator in the SPDE provides an alternative approach for
deriving stochastic Taylor expansions for the solution of such a SPDE.
Essentially, the exponential factor has a mollifying effect and ensures
that all integrals take values in the Hilbert space under
consideration. The iteration of such integrals allows us to derive
stochastic Taylor expansions of arbitrarily high order, which are
robust in the sense that they also hold for other types of driving
noise processes such as fractional Brownian motion. Combinatorial
concepts of trees and woods provide a compact formulation of the Taylor
expansions.
\end{abstract}

% KEYWORDS
%
\begin{keyword}[class=AMS]
\kwd{35K90}
\kwd{41A58}
\kwd{65C30}
\kwd{65M99}.
\end{keyword}
\begin{keyword}
\kwd{Taylor expansions}
\kwd{stochastic partial differential equations}
\kwd{SPDE}
\kwd{strong convergence}
\kwd{stochastic trees}.
\end{keyword}

\end{frontmatter}

%s1 ###
\section{Introduction}

Taylor expansions are a fundamental and repeatedly used method of
approximation in mathematics, in particular in numerical analysis.
Although numerical schemes for ordinary differential equations (ODEs)
are often derived in an ad hoc manner, those based on Taylor expansions
of the solution of an ODE, the Taylor schemes, provide a class of
schemes with known convergence orders against which other schemes can
be compared to determine their order. An important component of such
Taylor schemes are the iterated total derivatives of the vector field
corresponding higher derivatives of the solution, which are obtained
via the chain rule; see \cite{db}.

An analogous situation holds for It\^{o} stochastic ordinary
differential equations (SODEs), except, due to the less robust nature
of stochastic calculus, stochastic Taylor schemes instead of classical
Taylor schemes are the starting point to obtain consistent higher order
numerical schemes, see \cite{kp} for the general theory. Another
important difference is that SODEs are really just a symbolic
representation of stochastic integral equations since their solutions
are not differentiable, so an integral version of Taylor expansions
based on iterated application of the stochastic chain rule, the It\^{o}
formula, is required. Underlying this method is the fact that the
solution of a SODE is an It\^{o}-process or, more generally, a
semi-martingale and in particular of finite quadratic variation.

This approach fails, however, if a SODE is driven by an additive
stochastic process with infinite quadratic variation such as a
fractional Brownian motion, because the It\^{o} formula is in general
no longer valid. A new method to derive Taylor expansions in such cases
was presented in \cite{j,jk07,jk08a}. It uses the smoothness of the
coefficients, but only minimal assumptions on the nature of the driving
stochastic process. The resulting Taylor expansions there are thus
robust with respect to assumptions concerning the driving stochastic
process and, in particular, remain valid for other noise processes.

A similar situation holds for stochastic partial differential
equations\break
(SPDEs). In this article, we consider parabolic SPDEs with additive
noise of the form
%
%
%e1 ###
\begin{equation} \label{eqspde}
d U_t = [ A U_t + F( U_t ) ] \,dt + B \,d W_t ,\qquad U_0 = u_0,\qquad t \in[0,T],
\end{equation}
in a separable Hilbert space $H$, where $A$ is an unbounded linear
operator, $F$ is a nonlinear smooth function, $B$ is a bounded linear
operator and $W_t$, $ t \geq0 $, is an infinite-dimensional Wiener
process. (See Section \ref{sec2} for a precise description of the
equation above and the assumptions, we use.) Although the SPDE
(\ref{eqspde}) is driven by a martingale Brownian motion, the solution
process is with respect to a reasonable state space in general not a
semi-martingale anymore (see \cite{nt} for a clear discussion of this
problem) and except of special cases a general It\^{o} formula does not
exist for its solution (see, e.g., \cite{nt,pr}). Hence, stochastic
Taylor expansions for the solution of the SPDE (\ref{eqspde}) cannot be
derived as in \cite{kp} for the solutions of finite-dimensional SODEs.
Consequently, until recently, only low order numerical approximation
results for such SPDEs have been available (except for SPDEs with
spatially smooth noise; see, e.g., \cite{g03}).
%
% Consequently, until now, only temporal approximations
% of low order have been derived for the solutions of such SPDEs
%(except for SPDEs with finitely many
% stochastic processes, see e.g. \cite{g03}).
%
For example, the stochastic convolution of the semigroup generated by
the Laplacian with Dirichlet boundary conditions on the one-dimensional
domain $(0,1)$ and a cylindrical $I$-Wiener process on $ H = L^2 (
(0,1), \mathbb{R} ) $ has sample paths which are only $( \frac{1}{4} -
\varepsilon) $-H\"{o}lder continuous (see Section \ref{spacetime})
and previously considered approximations such as the linear implicit
Euler scheme or the linear implicit Crank--Nicolson scheme are not of
higher temporal order. The reason is that the infinite-dimensional
noise process has only minimal spatial regularity and the convolution
of the semigroup and the noise is only as smooth in time as in space.
This comparable regularity in time and space is a fundamental property
of the dynamics of semigroups; see \cite{sy} or also \cite{rn}, for
example. To overcome these problems, we thus need to derive robust
Taylor expansions for a SPDE of the form (\ref{eqspde}) driven by an
infinite-dimensional Brownian motion.

An idea used in \cite{jk08b} to derive what was called the exponential
Euler scheme for the SPDE (\ref{eqspde}), that has a better convergence
rate than hitherto analyzed schemes, can be exploited here. It is based
on the fact that the SPDE (\ref{eqspde}) can be understood in the mild
sense, that is as an integral equation of the form
%
%
%e2 ###
\begin{equation} \label{eqmain}
U_t = e^{ A t } u_0 + \int^{ t }_{ 0 } e^{ A (t-s) } F( U_s ) \,ds +
\int^{ t }_{ 0 } e^{ A (t-s) } B \,dW_s
\end{equation}
for all $ t \in[0,T] $ rather than as an integral equation obtained by
directly integrating the terms of the SPDE (\ref{eqspde}). (This mild
integral equation form of the SPDE is considered in some detail in the
monograph \cite{dpz}, (7.1) and (7.3.4), and in the monograph
\cite{pr}, (F.0.2).) The crucial point here is that all integrals in
the mild integral equation (\ref{eqmain}) contain the exponential
factor $e^{A (t-s)}$ of the operator $A$, which acts in a sense as a
mollifier and ensures that iterated versions of the terms remain in the
Hilbert space $H$. The main idea of the Taylor expansions presented in
this article is to use a classical Taylor expansion for the
nonlinearity $ F $ in the mild integral equation above and then to
replace the higher order terms recursively by Taylor expansions of
lower orders (see Section \ref{sec3}). Hence, this method avoids the need of an
It\^{o} formula but nevertheless yields stochastic Taylor expansions of
arbitrarily high order for the solution of the SPDE (\ref{eqspde}) (see
Section \ref{secabsex} for details). Moreover, these Taylor
expansions are robust with respect to the type of noise used and can
easily be modified to other types of noise such as fractional Brownian
motion.
% Finally, we note that in the case of finite dimensional
% SODEs driven by finite dimensional Brownian
% motions the resulting Taylor approximations obtained
% in this article are often identical to the
% approximations obtained via stochastic calculus and the It\^{o}
%formula (see Subsection \ref{SODEs}).

The paper is organized as follows. In the next section, we describe
precisely the SPDE that we are considering and state the assumptions
that we require on its terms and coefficients and on the initial value.
Then, in the third section, we sketch the idea and the notation for
deriving simple Taylor expansions, which we develop in section four in
some detail using combinatorial objects, specifically stochastic trees
and woods, to derive Taylor expansions of an arbitrarily high order. We
also provide an estimate for the remainder terms of the Taylor
expansions there. (Proofs are postponed to the final section.) These
results are illustrated with some representative examples in the fifth
section. Numerical schemes based on these Taylor expansions are
presented in the sixth section.
%will be discussed elsewhere.

%s2 ###
\section{Assumptions} \label{sec2}

Fix $ T > 0 $ and let $ ( \Omega, \mathcal{F}, \mathbb{P} ) $ be a
probability space with a normal filtration $ \mathcal{F}_t $, $ t
\in[0,T]$; see, for example, \cite{dpz} for details. In addition, let $ ( H, \langle\cdot,
\cdot\rangle ) $ be a separable $\mathbb{R}$-Hilbert space with its
norm denoted by $ | \cdot| $ and consider the SPDE (\ref{eqspde}) in
the mild integral equation form (\ref{eqmain}) on $H$, where $ W_t $, $
t \in[0,T]$, is a cylindrical $Q$-Wiener process with $ Q = I $ on
another separable $\mathbb{R}$-Hilbert space $ ( U, \langle\cdot, \cdot\rangle ) $
(space--time white noise), $ B \dvtx U \rightarrow H $ is a bounded
linear operator and the objects $ A $, $ F $ and $ u_0 $ are specified
through the following assumptions.
\begin{assumption}[(Linear operator $A$)]\label{A1}
Let $ I $ be a countable or finite set, let $ ( \lambda_i )_{ i
\in\mathcal{I} } \subset\mathbb{R}$ be a family of real numbers with $
\inf_{i \in\mathcal{I} } \lambda_i > - \infty$ and let $ ( e_i )_{ i
\in\mathcal{I} } \subset H $ be an orthonormal basis of $H$. Assume
that the linear operator $ A \dvtx D(A) \subset H \rightarrow H $ is
given by
\[
A v = \sum_{ i \in\mathcal{I} } -\lambda_i \langle e_i, v \rangle e_i
\]
for all $ v \in D(A) $ with $ D(A) = \{ v \in H | \sum_{ i
\in\mathcal{I} } | \lambda_i |^2 | \langle e_i, v \rangle |^2 < \infty \} $.
\end{assumption}
\begin{assumption}[(Drift term $F$)]\label{A2}
The nonlinearity $ F \dvtx H \rightarrow H $ is smooth and regular in
the sense that $ F $ is infinitely often Fr\'{e}chet differentiable and
its derivatives satisfy $ {\sup_{ v \in H } }| F^{(i)}(v) | < \infty $
for all $ i \in\mathbb{N} := \{ 1, 2, \ldots \} $.
\end{assumption}

Fix $ \kappa\geq0 $ with $ \sup_{ i \in\mathcal{I} } ( \kappa+
\lambda_i ) > 0 $ and let $ D ( ( \kappa- A )^r ) $, $ r \in\mathbb{R}
$, denote the domains of powers of the operator $ \kappa- A \dvtx D (
\kappa- A ) = D(A) \subset H \rightarrow H $, see, for example,
\cite{sy}.
\begin{assumption}[(Stochastic convolution)]\label{A3}
There exist two real numbers $ \gamma\in(0,1) $,
$\delta\in(0,\frac{1}{2}] $ and a constant $ C > 0 $ such that
\[
\int^T_0 | (\kappa-A)^{\gamma} e^{ A s } B |_{\mathrm{HS}}^2 \,ds < \infty,\qquad
\int^t_0 | e^{ A s } B |_{\mathrm{HS}}^2 \,ds \leq C t^{ 2 \delta}
\]
holds for all $ t \in[0,1] $, where \mbox{$| \cdot|_{\mathrm{HS}}$} denotes the
Hilbert--Schmidt norm for Hilbert--Schmidt operators from $ U $ to $H$.
\end{assumption}
\begin{assumption}[(Initial value $ u_0 $)]\label{A4}
The $\mathcal{F}_0/\mathcal{B} (D ( ( \kappa- A )^{\gamma} ) )
$-measurable mapping $ u_0 \dvtx \Omega\rightarrow D ( ( \kappa- A
)^{\gamma} ) $ satisfies $ \mathbb{E} | ( \kappa- A )^{\gamma} u_0 |^p
< \infty $ for every $ p \in[1,\infty) $, where $ \gamma\in(0,1) $ is
given in Assumption \ref{A3}.
\end{assumption}

Similar assumptions are used in the literature on the approximation of
this kind of SPDEs (see, e.g., Assumption H1--H3 in \cite{h} or
see also \cite{mgr,mgrb,mgrw,jk08b}).
%These assumptions are standard in the literature on the approximation
%of this kind of SPDEs (see \cite{jk08b,mgr,mgrb,mgrw}).
This setup also includes trace class noise (see Section
\ref{traceclass}) and finite-dimensional SODEs with additive noise (see
Section \ref{SODEs}).

Henceforth, we fix $t_0 \in[0,T)$ and denote by $ \mathcal{P} $ the set
of all adapted stochastic processes
\[
X \dvtx \Omega\rightarrow C( [t_0,T] , H ) \qquad\mbox{with } {\sup_{ t_0 \leq
t \leq T }} | X_t |_{ L^p } < \infty\ \forall p \geq1,
\]
and with continuous sample paths, where $ | Z |_{ L^p } := ( \mathbb{E}
| Z |^p )^{1/p} $ is the $ L^p $-norm of a random variable $ Z
\dvtx \Omega\rightarrow H $. Under Assumptions \ref{A1} and \ref{A3},
it is well known that the stochastic convolution
\[
\int^t_{0} e^{ A(t-s) } B \,dW_s,\qquad t \in[t_0,T],
\]
has an (up to indistinguishability) unique version with continuous
sample paths (see Lemma \ref{lemsc} in Section \ref{secsc}). From
now on, we fix such a version of the stochastic convolution.
Hence, under Assumptions \ref{A1}--\ref{A4} it is well known that there
is a pathwise unique adapted stochastic process $ U \dvtx
\Omega\rightarrow C ( [0,T], H ) $ with continuous sample paths, which
satisfies (\ref{eqmain}) (see Theorems 7.4 and 7.6 in
\cite{dpz}). Even more, this process satisfies
%
%
%e3 ###
\begin{equation} \label{boundedmom}
{\sup_{ 0 \leq t \leq T }} | ( \kappa- A )^\gamma U_t |_{ L^p } < \infty
\end{equation}
for all $ p \geq1 $.

%s3 ###
\section{Taylor expansions}\label{sec3}

In this section, we present the notation and the basic idea behind the
derivation of Taylor expansions. We write
\[
\Delta U_t := U_t - U_{ t_0 },\qquad \Delta t := t - t_0
\]
for $t \in[t_0,T] \subset[0,T]$, thus $ \Delta U $ denotes the
stochastic process $ \Delta U_t$, $t \in[t_0,T]$, in $\mathcal{P}$.
First, we introduce some integral operators and an expression relating
them and then we show how they can be used to derive some simple Taylor
expansions.

%s3.1 ###
\subsection{Integral operators}
\label{secintegralop}
Let $ j \in\{ 0,1,2,{1^*} \}$, where the indices
$ \{ 0,1,2 \}$ will label expressions containing only a constant value or no
value of the SPDE solution, while ${1^*}$ will label a certain integral with
time dependent values of the SPDE solution in the integrand.
Specifically, we define the stochastic processes $ I^0_j \in\mathcal{P}
$ by\looseness=-1
\[
I^0_j(t) := \cases{ \displaystyle( e^{ A \Delta t } - I ) U_{t_0}, &\quad $j=0$, \vspace*{2pt}\cr
\displaystyle\int^{t}_{t_0} e^{ A (t-s) } F ( U_{t_0} ) \,ds, &\quad $j = 1$, \vspace*{2pt}\cr
\displaystyle\int^{t}_{t_0} e^{ A (t-s) } B \,dW_s, &\quad $j = 2$, \vspace*{2pt}\cr
\displaystyle\int^{t}_{t_0} e^{ A (t-s) } F( U_{s} ) \,ds, &\quad $j = {1^*}$,}
\]
for each $ t \in[t_0,T] $. Note that the stochastic process $
\int^t_{t_0} e^{ A(t-s) } B \,dW_s $ for $ t \in[t_0,T] $ given by
\[
\int^t_{t_0} e^{A(t-s)} B \,dW_s = \int^{t}_0 e^{ A(t-s) } B \,dW_s - e^{ A
(t-t_0) } \biggl( \int^{t_0}_0 e^{ A ( t_0 - s ) } B \,dW_s \biggr)
\]
for every $ t \in[t_0,T] $ is indeed in $ \mathcal{P} $. Given $ i
\in\mathbb{N} $ and $ j \in\{ 1,{1^*} \} $, we then define the
$i$-multilinear symmetric mapping $I^i_j \dvtx \mathcal{P}^i
% $:=$
% $\underbrace{\mathcal{P} \times\ldots
% \times\mathcal{P}}_{\text{$i$-times}}$
\rightarrow\mathcal{P}$ by
\[
I^i_j[ g_1, \ldots, g_i ](t) := \frac{ 1 }{ i! } \int^{t}_{t_0} e^{ A
(t-s) } F^{(i)}( U_{t_0} ) ( g_1(s), \ldots, g_i(s) ) \,ds,
\]
when $j = 1$ and by
\begin{eqnarray*}
&&I^i_j[ g_1, \ldots, g_i ](t)
\\
&&\qquad:= \int^{t}_{t_0} e^{ A (t-s) } \biggl( \int^1_0 F^{(i)}( U_{t_0} + r \Delta
U_s ) ( g_1(s), \ldots, g_i(s) ) \frac{ (1-r)^{(i-1)} }{ (i-1)! } \,dr \biggr)
\,ds,
\end{eqnarray*}
when $j = {1^*}$ for all $ t \in[t_0,T] $ and all $ g_1, \ldots, g_i
\in \mathcal{P} $. One immediately checks that the stochastic
processes\vspace*{1pt}
$ I^0_j \in\mathcal{P} $, $ j \in\{ 0,1,2,1^* \} $, and the mappings $
I^i_j \dvtx \mathcal{P}^i \rightarrow \mathcal{P} $, $ i \in\mathbb{N}
$, $ j \in\{ 1, 1^* \} $, are well defined.

The solution process $U$ of (\ref{eqmain}) obviously satisfies
%
%
%e4 ###
\begin{equation} \label{eqmaindif}
\Delta U_t = ( e^{ A \Delta t } - I ) U_{ t_0 } + \int^{ t }_{ t_0 }
e^{ A (t-s) } F( U_s ) \,ds + \int^{ t }_{ t_0 } e^{ A (t-s) } B \,dW_s
\end{equation}
or, in terms of the above integral operators,
\[
\Delta U_t = I^0_0(t) + I^0_{{1^*}}(t) + I^0_{2}(t)
\]
for every $ t \in[t_0,T] $, which we can write symbolically in the
space $ \mathcal{P}$ as
%
%
%e5 ###
\begin{equation} \label{eqmain3}
\Delta U = I^0_0 + I^0_{{1^*}} + I^0_{2}.
\end{equation}
The stochastic processes $ I^i_j[g_1, \ldots, g_i] $ for $g_1, \ldots,
g_i \in\mathcal{P}$, $j \in\{ 0, 1 , 2 \} $ and $ i \in\{ 0, 1, 2,
\ldots\} $ only depend on the solution at time $ t = t_0 $. These terms
are therefore useful approximations for the solution process $ U_t $, $
t \in[t_0,T]$. However, the stochastic processes $ I^i_{{1^*}}[g_1,
\ldots, g_i] $ for $ g_1, \ldots, g_i \in\mathcal{P} $ and $ i \in\{ 0,
1,2, \ldots\} $ depend on the solution path $ U_t $ with $ t \in[t_0,T]
$. Therefore, we need a further expansion for these processes. For
this, we will use the important formula
%
%
%e6 ###
\begin{eqnarray} \label{formel1}
I^0_{{1^*}} & = & I^0_1
+ I^1_{{1^*}}[ \Delta U ] \nonumber\\[-8pt]\\[-8pt]
&=& I^0_1
+ I^1_{{1^*}}[ I^0_0 ]
+ I^1_{{1^*}}[ I^0_{{1^*}} ]
+ I^1_{{1^*}}[ I^0_{2} ],\nonumber
\end{eqnarray}
which is an immediate consequence of integration by parts and
(\ref{eqmain3}), and the iterated formula
%
%
%e7 ###
\begin{eqnarray} \label{formel2}
I^i_{{1^*}}[ g_1, \ldots, g_i ] & = & I^i_1[ g_1, \ldots, g_i ] + I^{ (
i + 1 )}_{{1^*}}[ \Delta U , g_1, \ldots, g_i ]
\nonumber\\
& = & I^i_1[ g_1, \ldots, g_i ] + I^{ ( i + 1 )}_{{1^*}}[ I^0_0 , g_1,
\ldots, g_i ]\\
&&{} + I^{ ( i + 1 )}_{{1^*}}[ I^0_{{1^*}} ,
g_1, \ldots, g_i ] + I^{ ( i + 1 )}_{{1^*}}[ I^0_{2} , g_1, \ldots, g_i]\nonumber
\end{eqnarray}
for every $g_1, \ldots, g_i \in\mathcal{P}$ and every $i \in\mathbb{N}$
(see Lemma \ref{lemformel} for a proof of the equations above).

%s3.2 ###
\subsection{Derivation of simple Taylor expansions}
\label{secsimple}

To derive a further expansion of~(\ref{eqmain3}), we insert formula
(\ref{formel1}) to the stochastic process $I^0_{{1^*}}$, that is,
\[
I^0_{{1^*}} = I^0_1 + I^1_{{1^*}}[ I^0_0 ] + I^1_{{1^*}}[ I^0_{{1^*}} ]
+ I^1_{{1^*}}[ I^0_{2} ]
\]
into (\ref{eqmain3}) to obtain
\begin{eqnarray*}
\Delta U = I^0_0 + ( I^0_1 + I^1_{{1^*}}[ I^0_0 ] + I^1_{{1^*}}[
I^0_{{1^*}} ] + I^1_{{1^*}}[ I^0_{2} ] ) + I^0_{2} ,
\end{eqnarray*}
which can also be written as
%
%
%e8 ###
\begin{equation}\label{tay1rest}
\Delta U = I^0_0 + I^0_1 + I^0_2 + I^1_{{1^*}}[ I^0_0 ] + I^1_{{1^*}}[
I^0_{{1^*}} ] + I^1_{{1^*}}[ I^0_{2} ] .
\end{equation}
If we can show that the double integral terms $ I^1_{{1^*}}[ I^0_0 ] $,
$ I^1_{{1^*}}[ I^0_{1^*} ] $ and $ I^1_{{1^*}}[ I^0_2 ] $ are
sufficiently small (indeed, this will be done in the next section),
then we obtain the approximation
%
%
%e9 ###
\begin{equation} \label{tay1restb}
\Delta U \approx I^0_0 + I^0_1 + I^0_2,
\end{equation}
or, using the definition of the stochastic processes $ I^0_0 $, $ I^0_1
$ and $ I^0_2 $,
\[
\Delta U_t \approx ( e^{ A \Delta t } - I ) U_{t_0} + \int^{t}_{t_0}
e^{ A (t-s) } F( U_{t_0} ) \,ds + \int^{t}_{t_0} e^{ A (t-s) } B \,dW_s
\]
for $ t \in[t_0,T] $. Hence,
%
%
%e10 ###
\begin{equation} \label{tay1}
U_t \approx e^{ A \Delta t } U_{t_0} + \biggl( \int^{\Delta t}_0 e^{ A s } \,ds
\biggr) F( U_{t_0} ) + \int^{t}_{t_0} e^{ A (t-s) } B \,dW_s
\end{equation}
for $ t \in[t_0,T] $ is an approximation for the solution of SPDE
(\ref{eqspde}). Since the right-hand side of (\ref{tay1}) depends on
the solution only at time $ t_0 $, it is the first nontrivial Taylor
expansion of the solution of the SPDE (\ref{eqspde}). The remainder
terms $ I^1_{{1^*}}[ I^0_0 ] $, $ I^1_{{1^*}}[ I^0_{1^*} ] $ and $
I^1_{{1^*}}[ I^0_2 ] $ of this approximation can be estimated by
\[
| I^1_{{1^*}}[ I^0_0 ](t) + I^1_{{1^*}}[ I^0_{{1^*}} ](t) +
I^1_{{1^*}}[ I^0_{2} ](t) |_{ L^2 } \leq C ( \Delta t )^{ ( 1 + \min(
\gamma, \delta ) ) }
\]
for every $ t \in[t_0,T]$ with a constant $C > 0$ and where $
\gamma\in(0,1) $ and $ \delta\in(0,\frac{1}{2}] $ are given in
Assumption \ref{A3} (see Theorem \ref{mainthm} in the next section for
details).

We write $ Y_t = O ( ( \Delta t )^{ r } ) $ with $ r > 0 $ for a
stochastic process $ Y \in\mathcal{P} $, if $ | Y_t |_{ L^2 } \leq C (
\Delta t )^r $ holds for all $ t \in[t_0,T] $ with a constant $ C > 0
$. Therefore, we have
\begin{eqnarray*}
&&U_t - \biggl( e^{ A \Delta t } U_{t_0} + \biggl( \int^{\Delta t}_0 e^{ A s } \,ds \biggr)
F( U_{t_0} ) + \int^{t}_{t_0} e^{ A (t-s) } B \,dW_s \biggr)
\\
&&\qquad= O \bigl( ( \Delta t )^{ ( 1 + \min( \gamma, \delta) ) } \bigr)
\end{eqnarray*}
or
%
%
%e11 ###
\begin{eqnarray} \label{onestep}
U_t &=& e^{ A \Delta t } U_{t_0} + \biggl( \int^{ \Delta t }_{0} e^{ A s } \,ds \biggr)
F( U_{t_0} ) + \int^{t}_{t_0} e^{ A (t-s) } B \,dW_s
\nonumber\\[-8pt]\\[-8pt]
&&{} + O \bigl( ( \Delta t )^{ ( 1 + \min( \gamma, \delta) ) } \bigr)
.\nonumber
\end{eqnarray}
The approximation (\ref{tay1}) thus has order $1 + \min(\gamma,\delta)$
in the above strong sense. It plays an analogous role to the simplest
strong Taylor expansion in \cite{kp} on which the Euler--Maruyama
scheme for finite-dimensional SODEs is based and was in fact used in
\cite{jk08b} to derive the exponential Euler scheme for the SPDE
(\ref{eqspde}). Note that the Euler--Maruyama scheme in \cite{kp}
approximates the solution of an SODE with additive noise locally with
order $\frac{3}{2}$. Here, in the case of trace class noise, we will
have $ \gamma= \frac{1}{2} - \varepsilon$, $ \delta= \frac{1}{2} $ (see
Section \ref{traceclass}), and therefore the exponential Euler
scheme for the SPDE (\ref{eqspde}) in \cite{jk08b} also approximates
the solution locally with order $ \frac{3}{2} - \varepsilon$ (see
Section \ref{tcexpeuler}), while other schemes in use, in particular
the linear implicit Euler scheme or the Crank--Nicolson scheme,
approximate the solution with order $ \frac{1}{2} $ instead of order $
\frac{3}{2} $ as in the finite-dimensional case. Therefore, the Taylor
approximation (\ref{onestep}) attains the classical order of the Euler
approximation for finite-dimensional SODEs and in general the Taylor
expansion introduced above lead to numerical schemes for SPDEs, which
converge with a higher order than other schemes in use (see Section
\ref{secnumerics}).

%s3.3 ###
\subsection{Higher order Taylor expansions} \label{sechigher}

Further expansions of the remainder terms in a Taylor expansion give a
Taylor expansion of higher order. To illustrate this, we will expand
the terms $ I^1_{{1^*}}[ I^0_{0} ] $ and $ I^1_{{1^*}}[ I^0_{2} ] $ in
(\ref{tay1rest}). From (\ref{formel2}), we have
\[
I^1_{{1^*}}[ I^0_{0} ]  = I^1_{1}[ I^0_{0} ] + I^2_{{1^*}}[ I^0_0,
I^0_{0} ] + I^2_{{1^*}}[ I^0_{{1^*}}, I^0_{0} ] + I^2_{{1^*}}[ I^0_{2},
I^0_{0} ]
\]
and
\[
I^1_{{1^*}}[ I^0_{2} ]  = I^1_{1}[ I^0_{2} ] + I^2_{{1^*}}[ I^0_0,
I^0_{2} ] + I^2_{{1^*}}[ I^0_{{1^*}}, I^0_{2} ] + I^2_{{1^*}}[ I^0_{2},
I^0_{2} ],
\]
which we insert into (\ref{tay1rest}) to obtain
\[
\Delta U = I^0_0 + I^0_1 + I^0_2 + I^1_1[ I^0_0 ] + I^1_1[ I^0_2 ] + R,
\]
where the remainder term $ R $ is given by
\begin{eqnarray*}
R & = & I^1_{{1^*}}[ I^0_{{1^*}} ] + I^2_{{1^*}}[ I^0_0, I^0_{0} ] +
I^2_{{1^*}}[ I^0_{{1^*}}, I^0_{0} ] + I^2_{{1^*}}[ I^0_{2}, I^0_{0} ]
\\
& &{} + I^2_{{1^*}}[ I^0_0, I^0_{2} ] + I^2_{{1^*}}[ I^0_{{1^*}}, I^0_{2}
] + I^2_{{1^*}}[ I^0_{2}, I^0_{2} ] .
\end{eqnarray*}
From Theorem \ref{mainthm} in the next section, we will see $R = O ( (
\Delta t )^{ ( 1 + 2 \min( \gamma, \delta) ) } ) $. Thus, we have
\[
\Delta U = I^0_0 + I^0_1 + I^0_2 + I^1_2[ I^0_0 ] + I^1_2[ I^0_2 ] + O
\bigl( ( \Delta t )^{ ( 1 + 2 \min( \gamma, \delta) ) } \bigr) ,
\]
which can also be written as
\begin{eqnarray*}
U_t &=& e^{ A \Delta t } U_{t_0} + \biggl( \int^{\Delta t}_0 e^{ A s } \,ds \biggr)
F( U_{t_0} ) + \int^{t}_{t_0} e^{ A (t-s) } B
\,dW_s \\
&&{} + \int^t_{t_0} e^{ A ( t-s) } F'( U_{t_0} ) ( e^{ A \Delta s } - I )
U_{ t_0}
\,ds \\
&&{} + \int^t_{t_0} e^{ A ( t-s) } F'( U_{t_0} ) \int^s_{t_0} e^{ A (
s-r) } B \,dW_r \,ds + O \bigl( ( \Delta t )^{ ( 1 + 2 \min( \gamma, \delta) ) }
\bigr) .
\end{eqnarray*}
This approximation is of order $ 1 + 2 \min( \gamma, \delta) $.

By iterating this idea, we can construct further Taylor expansions. In
particular, we will show in the next section how a Taylor expansion
of arbitrarily high order can be achieved.

%s4 ###
\section{Systematic derivation of Taylor expansions
of arbitrarily high order}

The basic mechanism for deriving a Taylor expansion for the SPDE (\ref
{eqspde}) was explained in the previous section. We illustrate now how
Taylor expansions of arbitrarily high order can be derived and will
also estimate their remainder terms. For this, we will identify the
terms occurring in a Taylor expansions by combinatorial objects, that is,
trees. It is a standard tool in numerical analysis to describe higher
order terms in a Taylor expansion via rooted trees (see, e.g.,
\cite{jb} for ODEs and \cite{b00,r01,r04,r06} for SODEs). In
particular, we introduce a class of trees which is appropriate for our
situation and show how the trees relate to the desired Taylor
expansions.

%s4.1 ###
\subsection{Stochastic trees and woods}
We begin with the definition of the trees that we need, adapting the
standard notation of the trees used in the Taylor expansion of SODEs
(see, e.g., Definition 2.3.1 in \cite{r01} as well as
\cite{b00,r04,r06}).

Let $ N \in\mathbb{N} $ be a natural number and let
\[
\mathbf{t}'\dvtx \{ 2, \ldots, N \} \rightarrow\{ 1, \ldots, N-1
\},\qquad
\mathbf{t}''\dvtx \{ 1, \ldots, N \} \rightarrow\{ 0, 1, 2, {1^*} \}
\]
be two mappings with the property that $ \mathbf{t}'(j) < j $ for all $
j \in\{ 2, \ldots, N \} $. The pair of mappings $ \mathbf{t} = (
\mathbf{t}', \mathbf{t}'' ) $ is a {S-tree} (stochastic tree) of length
$ N = l(\mathbf{t})$ nodes.

Every {S-tree} can be represented as a graph, whose nodes are given by
the set $ \operatorname{nd}( \mathbf{t}) := \{ 1, \ldots, N \} $ and
whose arcs are described by the mapping $ \mathbf{t}' $ in the sense
that there is an edge from $ j $ to $ \mathbf{t}'(j) $ for every node $
j \in\{ 2, \ldots, N \} $. In view of a rooted tree, $ \tau' $ also
codifies the parent and child pairings and is therefore often referred
as son-farther mapping (see, e.g., Definition 2.1.5 in \cite{r01}). The
mapping $ \mathbf{t}'' $ is an additional labeling of the nodes with $
\mathbf{t}''(j) \in\{ 0, 1, 2, {1^*} \} $ indicating the type of node $
j $ for every $ j \in\operatorname {nd}(\mathbf{t}) $.
%
%
%f1 ###
\begin{figure}[b]

\includegraphics{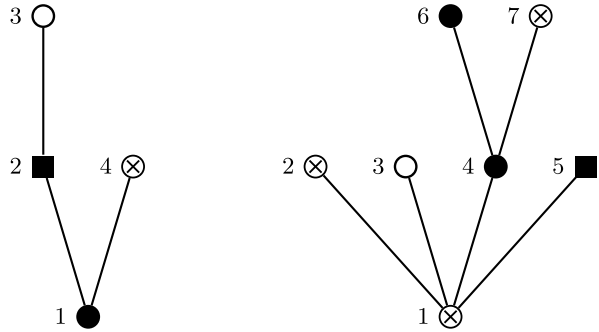}

\caption{Two examples of stochastic trees.}
\label{fig1}
\end{figure}
The left picture in Figure \ref{fig1} corresponds to the tree $
\mathbf{t}_1 = ( \mathbf{t}_1', \mathbf{t}_1'' ) $ with $
\operatorname{nd}( \mathbf{t}_1 ) = \{ 1,2,3,4 \} $ given by
\[
\mathbf{t}_1'(4) = 1,\qquad
\mathbf{t}_1'(3) = 2,\qquad
\mathbf{t}_1'(2)= 1
\]
and
\[
\mathbf{t}_1''(1) = 1,\qquad
\mathbf{t}_1''(2)= {1^*},\qquad
\mathbf{t}_1''(3)= 2,\qquad
\mathbf{t}_1''(4)= 0.
\]
The root is always presented as the lowest node. The number on the left
of a node in Figure \ref{fig1} is the number of the node of the
corresponding tree. The type of the nodes in Figure \ref{fig1} depends
on the additional labeling of the nodes given by~$ \mathbf{t}_1'' $.
More precisely, we represent a node $ j \in\operatorname{nd}(\mathbf
{t}_1) $ by
\ \psdot*[dotsize=9pt, dotstyle=otimes](0pt,3pt)\mbox{ }\
if $ \mathbf{t}_1''(j) = 0 $, by
\raisebox{0.6ex}[0pt][0pt]{\TC*[radius=4.5pt]}
if $ \mathbf {t}_1''(j) = 1 $, by
\raisebox{0.6ex}[0pt][0pt]{\TC[linewidth=0.9pt,radius=4.5pt]}
if $ \mathbf{t}_1''(j) = 2 $, and finally by
\raisebox{0.6ex}[0pt][0pt]{\Tf*[framesize=8pt]}
if $ \mathbf{t}_1''(j) = {1^*} $.
The right picture in Figure \ref{fig1} corresponds to the tree $
\mathbf{t}_2 = ( \mathbf{t}_2', \mathbf{t}_2'' ) $ with $
\operatorname{nd}( \mathbf{t}_2 ) = \{ 1, \ldots, 7 \} $ given by
\begin{eqnarray*}
\mathbf{t}_2'(7) &=& 4,\qquad
\mathbf{t}_2'(6) = 4,\qquad
\mathbf{t}_2'(5) = 1,\\
\mathbf{t}_2'(4) &=& 1,\qquad
\mathbf{t}_2'(3) = 1,\qquad
\mathbf{t}_2'(2)= 1
\end{eqnarray*}
and
\begin{eqnarray*}
\mathbf{t}_2''(1) &=& 0,\qquad
\mathbf{t}_2''(2) = 0,\qquad
\mathbf{t}_2''(3) = 2,\qquad
\mathbf{t}_2''(4) = 1,\\
\mathbf{t}_2''(5) &=& {1^*},\qquad
\mathbf{t}_2''(6)= 1,\qquad
\mathbf{t}_2''(7)= 0.
\end{eqnarray*}
We denote the set of all stochastic trees by $\mathbf{ST}$ and will
also consider a tuple of trees, that is, a wood. The set of {S-woods}
({stochastic woods}) is defined by
\[
\mathbf{SW} := \bigcup^{\infty}_{n=1} ( \mathbf{ST} )^n.
\]
Of course, we have the embedding $ \mathbf{ST} \subset\mathbf{SW} $. A
simple example of an {S-wood} which will be required later is $
\mathbf{w} _0 = (\mathbf{t}_1, \mathbf{t}_2, \mathbf{t}_3) $ with $
\mathbf{t}_1 $, $ \mathbf{t}_2 $ and $ \mathbf{t}_3 $ given by $
l(\mathbf{t}_1) = l(\mathbf{t}_2) = l(\mathbf {t}_3) = 1 $ and $
\mathbf{t}_1''(1) = 0 $, $ \mathbf{t}_2''(1) = {1^*} $, $
\mathbf{t}_3''(1) = 2 $. This is shown in Figure~\ref{fig4}
%
%f2 ###
\begin{figure}

\includegraphics{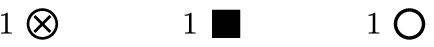}

\caption{The stochastic wood $\mathbf{w}_0 $ in $\mathbf{SW}$.}
\label{fig4}
\end{figure}
where the left tree corresponds to $ \mathbf{t}_1 $, the middle one to
$\mathbf{t}_2 $ and the right tree corresponds to $ \mathbf{t}_3 $.

%s4.2 ###
\subsection{Construction of stochastic trees and woods}

We define an operator on the set $\mathbf{SW}$, that will enable us to
construct an appropriate stochastic wood step by step. Let $ \mathbf{w}
= ( \mathbf{t}_1, \ldots, \mathbf{t}_n ) $ with $ n \in\mathbb{N}$ be a
{S-wood} with $ \mathbf{t}_i = ( \mathbf{t}_i', \mathbf{t}_i'' )
\in\mathbf{ST} $ for every $ i \in\{ 1, 2, \ldots, n \} $. Moreover,
let $ i \in\{ 1, \ldots, n \} $ and $ j \in\{ 1, \ldots,
l(\mathbf{t}_i) \} $ be given and suppose that $ \mathbf{t}_i''( j ) =
{1^*}$, in which case we call the pair $ (i,j) $ an active node of $
\mathbf{w}$. We denote the set of all active nodes of $ \mathbf{w}$ by
$\operatorname{acn}( \mathbf{w})$.

\begin{figure}

\includegraphics{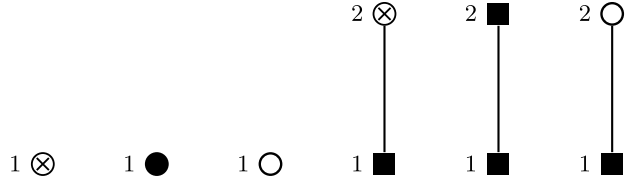}

\caption{The stochastic wood $\mathbf{w}_1 $ in $\mathbf{SW}$.}
\label{fig2}
\end{figure}

Now, we introduce the trees $ \mathbf{t}_{n+1} = ( \mathbf{t}_{n+1}',
\mathbf{t}_{n+1}'' ) $, $ \mathbf{t}_{n+2} = ( \mathbf{t}_{n+2}',
\mathbf{t}_{n+2}'' ) $ and $ \mathbf{t}_{n+3} = ( \mathbf{t}_{n+3}',
\mathbf{t}_{n+3}'' ) $ in $\mathbf{ST}$ by $ \operatorname{nd}(
\mathbf{t}_{n+m} ) = \{ 1, \ldots, l(\mathbf{t}_i), l(\mathbf{t}_i)+1
\} $ and
\begin{eqnarray*}
\mathbf{t}_{n+m}'( k ) &=& \mathbf{t}_i'( k ),\qquad
k=2,\ldots,l(\mathbf{t}_i),\hspace*{46pt}\\
\mathbf{t}_{n+m}''( k ) &=& \mathbf{t}_i''( k ),\qquad
k=1,\ldots,l( \mathbf{t}_i ) ,\hspace*{46pt}
\end{eqnarray*}
\begin{eqnarray*}
\mathbf{t}_{n+m}' \bigl( l(\mathbf{t}_i)+1 \bigr) &=& j,\qquad \mathbf{t}_{n+m}'' \bigl(
l(\mathbf{t}_i)+1 \bigr) = \cases{{1^*}, &\quad $m=2$, \cr (m-1),
&\quad else,}\hspace*{-27pt}\hspace*{10pt}
\end{eqnarray*}
for $ m = 1,2,3 $. Finally, we consider the {S-tree} $
\tilde{\mathbf{t}}^{(i,j)} = ( \tilde{\mathbf{t}}', \tilde
{\mathbf{t}}'' ) $ given by $ \tilde{ \mathbf{t}}' = \mathbf{t}_i' $,
but with $ \tilde{ \mathbf{t}}''(k) = \mathbf{t}_i''(k) $ for $ k \neq
j $ and $ \tilde{ \mathbf{t}}''(j) = 1 $. Then, we define
\[
E_{(i,j)}\mathbf{w} = E_{(i,j)} ( \mathbf{t}_1, \ldots, \mathbf{t}_n )
:= \bigl( \mathbf{t}_1, \ldots, \mathbf{t}_{i-1},
\tilde{\mathbf{t}}^{(i,j)}, \mathbf{t}_{i+1}, \ldots, \mathbf{t}_{n+3}
\bigr)
\]
and consider the set of all woods that can be constructed by
iteratively applying the $ E_{(i,j)} $ operations, that is, we define
\begin{eqnarray*}\hspace*{2pt}
\mathbf{SW}' &:=& \{ \mathbf{w}_0 \}
\cup \left\{\hspace*{-1pt} \mathbf{w}\in\mathbf{SW} \left|
\begin{array}{l}
\exists n \in\mathbb{N}, i_1, j_1, \ldots, i_n, j_n
\in\mathbb{N} \dvtx \\
\forall l=1,\ldots,n \dvtx ( i_l, j_l ) \in\operatorname{acn}\bigl( E_{
(i_{l-1}, j_{l-1}) } \cdots E_{(i_1,j_1)} \mathbf{w}_0 \bigr) ,
\\
\hspace*{2.1pt}\mathbf{w}= E_{(i_n,j_n)} \cdots E_{(i_1,j_1)} \mathbf{w}_0
\end{array}\hspace*{-5pt}
\right.\right\}
\end{eqnarray*}
for the $\mathbf{w}_0$ introduced above. To illustrate these
definitions, we present some examples using the initial stochastic wood
$ \mathbf{w}_0 $ given in Figure \ref{fig4}. We present these examples
here in a brief way, and, later in Section \ref{secabsex}, we
describe more detailed the main advantages of the particular examples
considered here. First, the active nodes of $ \mathbf{w}_0 $ are
$\operatorname{acn}( \mathbf{w}_0 ) = \{ (2,1) \}$, since the first
node in the second tree in $ \mathbf{w}_0 $ is the only node of type $
{1^*} $. Hence, $ E_{(2,1)} \mathbf{w}_0$ is well defined and the
resulting stochastic wood $ \mathbf{w}_1 = E_{(2,1)} \mathbf{w}_0 $,
which has six trees, is presented in Figure \ref {fig2}. Writing $
\mathbf{w}_1 = ( \mathbf{t}_1, \ldots, \mathbf {t}_6 )$, the left tree
in Figure \ref{fig2} corresponds to $ \mathbf{t}_1 $, the second tree
in Figure \ref{fig2} corresponds to $ \mathbf{t}_2 $, and so on.
%
%
%f3 ###
%
%
%f4 ###
\begin{figure}[b]

\includegraphics{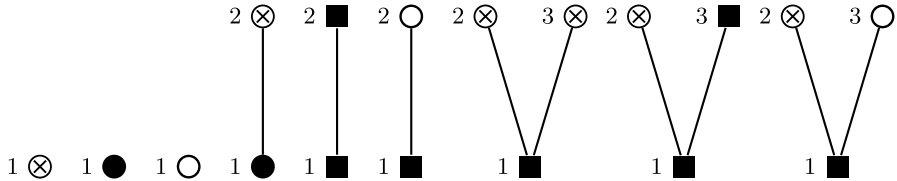}

\caption{The stochastic wood $\mathbf{w}_2$ in $\mathbf{SW}$.}
\label{fig5}
\end{figure}
Moreover, we have
%
%
%e12 ###
\begin{equation} \label{anodefig2}
\operatorname{acn}( \mathbf{w}_1 ) = \{ (4,1), (5,1), (5,2), (6,1) \}
\end{equation}
for the active nodes of $ \mathbf{w}_1 $, so $ \mathbf{w}_2 = E_{(4,1)}
\mathbf{w}_1 $ is also well defined. It is presented in Figure~\ref{fig5}.
In Figure \ref{fig6}, we present the stochastic wood
$ \mathbf{w}_3 = E_{ (6,1) } \mathbf{w}_2 $,
%
%
%f5 ###
\begin{figure}

\includegraphics{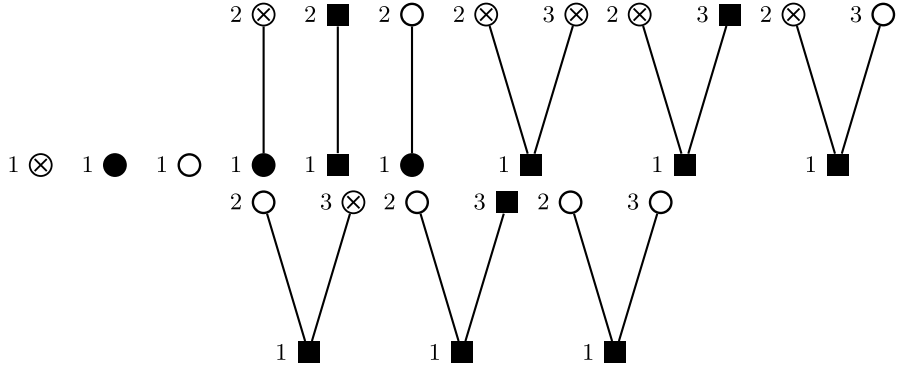}

\caption{The stochastic wood $\mathbf{w}_3 $
in $\mathbf{SW}$.}
\label{fig6}
\end{figure}
which is well defined since
%
%
%e13 ###
\begin{equation} \label{anodefig5}
\operatorname{acn}( \mathbf{w}_2 ) = \{ (5,1), (5,2), (6,1), (7,1),
(8,1), (8,3), (9,1) \} .
\end{equation}
For the S-wood $ \mathbf{w}_3 $, we have
%
%
%e14 ###
\begin{equation} \label{anodefig6}
\operatorname{acn}( \mathbf{w}_3 ) = \left\{
\matrix{
(5,1), (5,2), (7,1), (8,1), (8,3), (9,1), \cr
(10,1), (11,1), (11,3), (12,1)}
\right\} .
\end{equation}
Since $ (7,1) \in\operatorname{acn}(\mathbf{w}_3) $, the stochastic
wood $ \mathbf{w}_4 = E_{ (7,1) } \mathbf{w}_3 $ is well defined and
presented in Figure \ref{fig7}.
%
%
%f6 ###
\begin{figure}[b]

\includegraphics{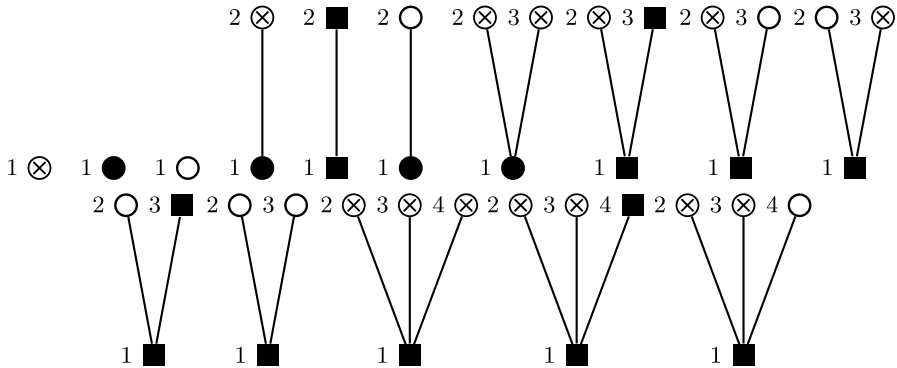}

\caption{The stochastic wood $\mathbf{w}_4 $
in $\mathbf{SW}$.}
\label{fig7}
\end{figure}
For the active nodes, we obtain
%
%
%e15 ###
\begin{equation} \label{anodefig7}
\operatorname{acn}( \mathbf{w}_4 ) =
\left\{
\matrix{
(5,1), (5,2), (8,1), (8,3),
(9,1), (10,1), (11,1), \cr
(11,3), (12,1),
(13,1), (14,1), (14,4), (15,1)
}
\right\} .
\end{equation}
Finally, we present the stochastic wood $ \mathbf{w}_5 = E_{(12,1)}
E_{(10,1)} E_{(9,1)} \mathbf{w}_4 $ with
%
%
%e16 ###
\begin{equation} \label{anodefig8}
\operatorname{acn}( \mathbf{w}_5 ) =
\left\{
\matrix{
(5,1), (5,2), (8,1), (8,3), (11,1),
(11,3), (13,1), (14,1), \cr
(14,4), (15,1), (16,1),
(17,1), (17,4), (18,1),
(19,1), \cr
(20,1), (20,4), (21,1),
(22,1), (23,1), (23,4), (24,1)}
\right\}\hspace*{-30pt}
\end{equation}
in Figure \ref{fig8}.
%
%
%f7 ###
\begin{figure}

\includegraphics{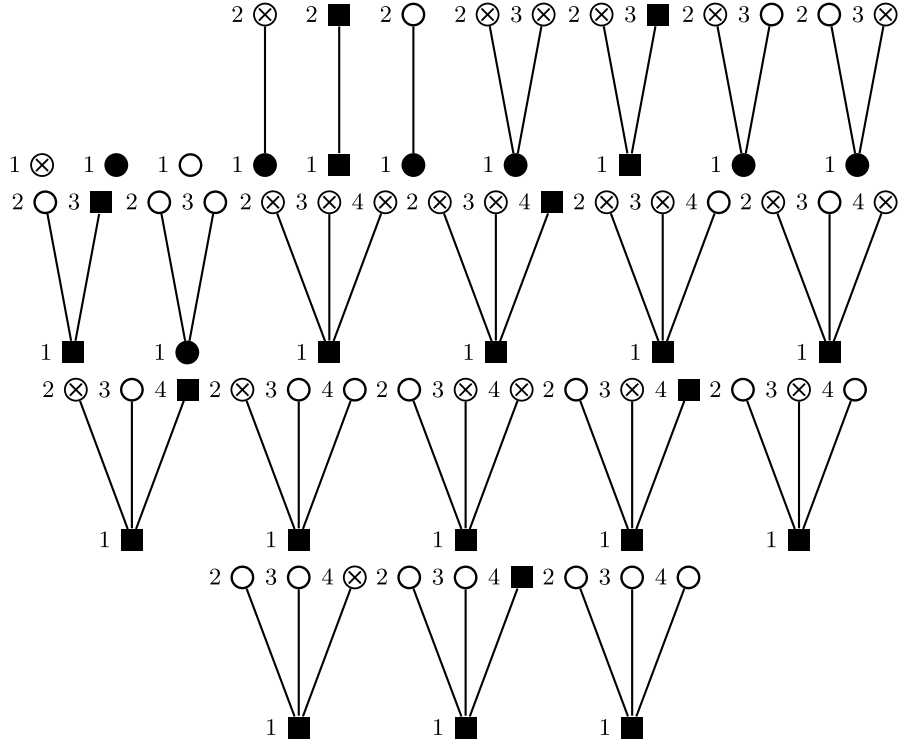}

\caption{The stochastic wood $\mathbf{w}_5$
in $\mathbf{SW}$.}
\label{fig8}
\end{figure}
By definition, the S-woods $\mathbf{w}_0, \mathbf{w}_1, \ldots, \mathbf{w}_5 $ are
in $\mathbf{SW}'$, but the stochastic wood given in Figure \ref{fig1}
is not in $\mathbf{SW}'$.

%s4.3 ###
\subsection{Subtrees}

Let $ \mathbf{t}= (\mathbf{t}',\mathbf{t}'') $ be a given S-tree with $
l( \mathbf{t}) \geq2 $. For two nodes $ k, l
\in\operatorname{nd}(\mathbf{t}) $ with $ k \leq l $, we say that $ l $
is a grandchild of $ k $ if there exists a sequence $ k_1 = k < k_2 <
\cdots< k_n = l$ of nodes with $ n \in\mathbb{N} $ such that $
\mathbf{t}' ( k_{v+1} ) = k_v $ for every $ v \in\{ 1, \ldots, n-1 \}
$. Suppose now that $ j_1, \ldots, j_n \in\operatorname{nd}(\mathbf{t})
$ with $ n \in\mathbb{N} $ and $ j_1 < \cdots< j_n $ are the nodes of $
\mathbf{t}$ such that $ \mathbf{t}'( j_i ) = 1 $ for every $i = 1,
\ldots, n$. Moreover, for a given $i \in \{ 1, \ldots, n \}$ suppose
that $ j_{i,1}, j_{i,2}, \ldots, j_{i,l_i} \in\operatorname{nd}(
\mathbf{t}) $ with $ j_i = j_{i,1} < j_{i,2} < \cdots< j_{i,l_i} \leq
l(\mathbf{t}) $ and $ l_i \in\mathbb{N}$ are the grandchildren of
$j_i$. Then, we define a tree $ \mathbf{t}_i = ( \mathbf{t}_i',
\mathbf{t}_i'' ) \in\mathbf{ST} $ with $ l( \mathbf{t}_i ) := l_i $ and
\[
j_{ i, \mathbf{t}_i'( k ) } = \mathbf{t}'( j_{i,k} ),\qquad
\mathbf{t}_i''( k) = \mathbf{t}''( j_{i,k} )
\]
for all $ k \in\{ 2, \ldots, l_i \} $ and $ \mathbf{t}_i''( 1 ) =
\mathbf{t}''( j_i ) $. We call the trees $ \mathbf{t}_1, \ldots,
\mathbf{t}_n \in\mathbf{ST} $ defined in this way the
\textit{subtrees} of $ \mathbf{t}$. For example, the subtrees of the
%
%
%f8 ###
\begin{figure}

\includegraphics{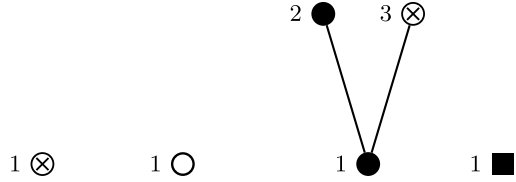}

\caption{Subtrees of the right tree in Figure \protect\ref{fig1}.}
\label{fig3}
\end{figure}
right tree in Figure \ref{fig1} are presented in Figure \ref{fig3}.

\eject
%s4.4 ###
\subsection{Order of a tree} \label{secorder}

Later stochastic woods in $\mathbf{SW}'$ will represent Taylor
expansions and Taylor approximations of the solution process $U$ of the
SPDE~(\ref{eqspde}).
% To this end stochastic woods in
% \textbf{ST} will represent certain
% summands within the Taylor
% approximations and expansions.
Additionally, we will estimate the approximation orders of these Taylor
approximations. To this end, we introduce the order of a stochastic
tree and of a stochastic wood, which is motivated by Lemma
\ref{importantlem} below. More precisely, let $\operatorname{ord} \dvtx
\mathbf{ST} \rightarrow[0,\infty )$ be given by
\begin{eqnarray*}
\operatorname{ord}( \mathbf{t}) & := & l( \mathbf{t}) + (\gamma- 1) |
\{ j \in\operatorname{nd}( \mathbf {t}) | \mathbf{t}''(j)
= 0 \} | \\
& &{} + (\delta- 1) | \{ j \in\operatorname{nd}( \mathbf{t}) | \mathbf{t}
''(j) = 2 \} |
\end{eqnarray*}
for every S-tree $ \mathbf{t}= ( \mathbf{t}', \mathbf{t}'')
\in\mathbf{ST} $. For example, the order of the left tree in
Figure~\ref{fig1} is $ 2 + \gamma+ \delta$ and the order of the right tree in
Figure \ref{fig1} is $ 3 + 3 \gamma+ \delta$ (since the right tree has
three  nodes of type $ 0 $, three nodes of type $ 1 $, respectively, $
{1^*} $, and one node of type $ 2 $).

In addition, we say that a tree $ \mathbf{t}= ( \mathbf{t}', \mathbf
{t}'') $ in $\mathbf{ST}$ is \textit{active} if there is a $ j
\in\operatorname{nd}( \mathbf{t}) $ such that $ \mathbf{t}''(j) = {1^*}
$. In that sense a S-tree is active if it has an active node. Moreover,
we define the order of an S-wood $ \mathbf{w} = ( \mathbf{t}_1, \ldots,
\mathbf{t}_n ) \in \mathbf{SW} $ with $ n \in\mathbb{N} $ as
\[
\operatorname{ord}( \mathbf{w}) := \min \{ \operatorname{ord}(
\mathbf{t}_i ), 1 \leq i \leq n | \mathbf{t}_i \mbox{ is active} \} .
\]
To illustrate this definition, we calculate the order of some
stochastic woods. First of all, the stochastic wood in Figure
\ref{fig4} has order $ 1 $, since only the middle tree in Figure
\ref{fig4} is active. More precisely, the node $ (2,1) $ of the S-wood
$ \mathbf{w}_0 $ is an active node and therefore the second tree is
active. The second tree in Figure \ref{fig4} has order~$ 1 $ (since it
only consists of one node of type $ 1^* $). Hence, the S-wood $
\mathbf{w}_0 $ has order~$ 1 $. Since the last three trees are active
in the stochastic wood $ \mathbf{w}_1 $ in Figure \ref{fig2} [see
(\ref{anodefig2}) for the active nodes of $\mathbf{w}_1 $], we obtain
that the stochastic wood in Figure \ref{fig2} has order $ 1 +
\min(\gamma,\delta) $. The last three trees in the S-wood $
\mathbf{w}_1 $ have order $ 1 + \gamma$, $ 2 $ and $ 1 + \delta$,
respectively. As a third example, we consider the S-wood $ \mathbf{w}_2
$ in Figure \ref{fig5}. The active nodes of $ \mathbf{w}_2 $ are
presented in (\ref{anodefig5}). Hence, the last five S-trees are
active. They have the orders $ 2 $, $ 1 + \delta$, $ 1 + 2 \gamma$, $ 2
+ \gamma$ and $ 1 + \gamma+ \delta$. The minimum of the five real
numbers is $ 1 + \min( 2 \gamma, \delta) $. Therefore, the order of the
S-wood $ \mathbf{w}_2 $ in Figure \ref{fig5} is $ 1 + \min( 2 \gamma,
\delta) $. A similar calculation shows that the order of the stochastic
wood $ \mathbf{w}_3 $ in Figure \ref{fig6} is $ 1 + 2
\min(\gamma,\delta) $ and that the order of the stochastic wood $
\mathbf{w}_4 $ in Figure \ref{fig7} is $ 1 + \min( 3 \gamma, \gamma+
\delta, 2 \delta) $. Finally, we obtain that the stochastic wood $
\mathbf{w}_5 $ in Figure \ref{fig8} is of order $ 1 + 3 \min(\gamma,
\delta,\frac{1}{3} ) $.

%s4.5 ###
\subsection{Trees and stochastic processes}

To identify each tree in $\mathbf{ST}$ with a predictable stochastic
process in $\mathcal{P}$, we define two functions $\Phi\dvtx
\mathbf{ST} \rightarrow\mathcal{P}$ and $\Psi\dvtx \mathbf{ST}
\rightarrow\mathcal{P}$, recursively. For a given S-tree $\mathbf{t} =
(\mathbf{t}',\mathbf{t}'') \in\mathbf{ST}$, we define $\Phi(
\mathbf{t}) := I^0_{ \mathbf{t}''( 1 ) } $ when $ \mathbf{t}''( 1 ) \in
\{ 0,2 \} $ or $l( \mathbf{t}) = 1 $ and, when $ l( \mathbf{t}) \geq2$
and $\mathbf{t}''( 1 ) \in \{ 1, 1^* \}$, we define
\[
\Phi( \mathbf{t}) := I^n_{ \mathbf{t}''( 1 ) } [ \Phi( \mathbf{t}_1 ),
\ldots, \Phi( \mathbf{t}_n ) ],
\]
where $\mathbf{t}_1, \ldots, \mathbf{t}_n \in\mathbf{ST} $ with $n
\in\mathbb{N}$ are the subtrees of $\mathbf{t}$. In addition, for an
arbitrary $ \mathbf{t} \in\mathbf{ST}$, we define $\Psi( \mathbf{t}) :=
0$ if $\mathbf{t}$ is an active tree and $\Psi( \mathbf{t}) = \Phi(
\mathbf{t})$ otherwise. Finally, for a S-wood $ \mathbf{w} = (
\mathbf{t}_1, \ldots, \mathbf{t}_n )$ with $n \in\mathbb{N}$ we define
$ \Phi( \mathbf{w})$ and $ \Psi( \mathbf{w}) $ by
\[
\Phi( \mathbf{w}) = \Phi( \mathbf{t}_1 ) + \cdots+ \Phi( \mathbf {t}_n
),\qquad \Psi( \mathbf{w}) = \Psi( \mathbf{t}_1 ) + \cdots+ \Psi( \mathbf
{t}_n ) .
\]

As an example, we have
%
%
%e17 ###
\begin{equation} \label{extr1a}
\Phi( \mathbf{w}_0 ) = I^0_0 + I^0_{{1^*}} + I^0_2 \quad\mbox{and}\quad \Psi(
\mathbf{w}_0 ) = I^0_0 + I^0_2
\end{equation}
for the elementary stochastic wood $ \mathbf{w}_0 $ (see Figure
\ref{fig4}). Hence, we obtain
%
%
%e18 ###
\begin{equation} \label{treesol1}
\Phi( \mathbf{w}_0 ) = \Delta U
\end{equation}
from (\ref{eqmain3}) and (\ref{extr1a}). Since $(2,1)$ is an active
node of $ \mathbf{w}_0 $, we obtain
%
%
%e19 ###
\begin{equation} \label{extr1}
\Phi( \mathbf{w}_1 ) = I^0_0 + I^0_{1} + I^0_2 + I^1_{{1^*}}[ I^0_0 ] +
I^1_{{1^*}}[ I^0_{{1^*}} ] + I^1_{{1^*}}[ I^0_2 ]
\end{equation}
and
%
%
%e20 ###
\begin{equation} \label{extr2}
\Psi( \mathbf{w}_1 ) = I^0_0 + I^0_{1} + I^0_2
\end{equation}
for the S-wood $ \mathbf{w}_1 = E_{(2,1)} \mathbf{w}_0 $ presented in
Figure \ref{fig2}. Moreover, in view of (\ref{formel1}) and~(\ref{formel2}), we have
%
%
%e21 ###
\begin{equation} \label{eformel}
\Phi( \mathbf{w}) = \Phi\bigl( E_{(i,j)} \mathbf{w}\bigr)
\end{equation}
for every active node $ (i,j) \in\operatorname{acn}(\mathbf{w}) $ and
every stochastic wood $\mathbf{w} \in\mathbf{SW}'$.

Hence, we obtain
\[
\Phi( \mathbf{w}_1 ) = \Phi\bigl( E_{(2,1)} \mathbf{w}_0 \bigr)
= \Phi(\mathbf{w}_0 )
\]
due to the equation above and the definition of $ \mathbf{w}_1 $.
Hence, we obtain $ \Phi( \mathbf{w}_1 ) = \Delta U, $ which can also be
seen from (\ref{extr1}), since the right-hand side of (\ref{extr1}) is
nothing other than (\ref{tay1rest}). We also note that the right-hand
side of (\ref{extr2}) is just the exponential Euler approximation in
(\ref{tay1restb}), so we obtain
\[
\Delta U = \Phi( \mathbf{w}_1 ) \approx\Psi( \mathbf{w}_1 ) .
\]
With the above notation and definitions we are now able to present the
main result of this article, which is a representation formula for the
solution of the SPDE (\ref{eqspde}) via Taylor expansions and an
estimate of the remainder terms occurring in the Taylor expansions.
\begin{theorem} \label{mainthm}
Let Assumptions \ref{A1}--\ref{A4} be fulfilled and let $\mathbf{w}\in
\mathbf{SW}'$ be an arbitrary stochastic wood. Then, for each $p
\in[1,\infty)$, there is a constant $ C_p > 0 $ such that
%
%
%e22 ###
\begin{eqnarray}\label{taylorSPDE}
U_t &=& U_{t_0} + \Phi(\mathbf{w})(t),\nonumber\\[-8pt]\\[-8pt]
\bigl( \mathbb{E} [ | U_t - U_{t_0} -
\Psi(\mathbf{w})(t) |^p ] \bigr)^{{1/p}}
&\leq&
C_p \cdot ( t - t_0)^{ \operatorname{ord}(\mathbf{w})}\nonumber
\end{eqnarray}
holds for every $ t \in[t_0,T] $, where $ U_t $, $ t \in[0,T] $, is the
solution of the SPDE (\ref{eqspde}). Here the constant $ C_p > 0 $ is
independent of $t$ and $t_0$ but depends on $ p$ as well as
$\mathbf{w}$, $T$ and the coefficients of the SPDE (\ref{eqspde}).
\end{theorem}

The representation of the solution here is a direct consequence of
(\ref{treesol1}) and (\ref{eformel}). The proof for the estimate in
(\ref{taylorSPDE}) will be given in Section \ref{secproofs}. Here, $
\Phi( \mathbf{w}) = \Delta U $ is the increment of the solution of the
SPDE (\ref{eqspde}), while $ \Psi( \mathbf{w}) $ is the Taylor
approximation of the increment of the solution and $\Phi( \mathbf{w}) -
\Psi( \mathbf{w})$ is its remainder for every arbitrary $
\mathbf{w}\in\mathbf{SW}' $. Since there are woods in $\mathbf{SW}'$
with arbitrarily high orders, Taylor expansions of arbitrarily high
order can be constructed by successively applying the operator $
E_{(i,j)}$ to the initial S-wood $ \mathbf{w}_0 $. Finally, the
approximation result of Theorem \ref{mainthm} can also be written as
\[
U_t = U_{ t_0 } + \Psi( \mathbf{w}) + O \bigl( ( \Delta t )^{
\operatorname{ord}(\mathbf{w}) } \bigr)
\]
for every stochastic wood $
\mathbf{w}\in\mathbf{SW}' $. Here, we also remark that Assumptions
\ref{A1}--\ref{A4} can be weakened. In particular, instead of
Assumption \ref{A3}, one can assume that the nonlinearity $ F \dvtx V
\rightarrow V $ is only $ i $-times Fr\'{e}chet differentiable with $ i
\in\mathbb{N} $ sufficiently high and that the derivatives of $ F $
satisfy only local estimates, where $ V \subset H $ is a continuously
embedded Banach space. Nevertheless, it is usual to present Taylor
expansions for stochastic differential equations under such restrictive
assumptions here (see \cite{kp}) and then after considering a
particular numerical scheme one reduces these assumptions by pathwise
localization techniques (see \cite{g98b,jkn} for SDEs and \cite{j09b}
for SPDEs).

%s5 ###
\section{Examples} \label{secex}

We present some examples here to illustrate the Taylor expansions
introduced above.

%s5.1 ###
\subsection{Abstract examples of the Taylor expansions} \label{secabsex}
We begin with some abstract examples of the Taylor expansions.

%s5.1.1 ###
\subsubsection{Taylor expansion of order $1$} \label{sec511}
The first Taylor expansion of the solution is given by the initial
stochastic wood $ \mathbf{w}_0 $ (see Figure \ref{fig4}), that is, we have
$ \Phi( \mathbf{w}_0 ) = \Delta U $ approximated by $ \Psi(
\mathbf{w}_0 ) $ with order $ \operatorname{ord}( \mathbf {w}_0 ) $.
Precisely, we have
\[
\Psi( \mathbf{w}_0 )(t) = ( e^{ A \Delta t } - I ) U_{ t_0 } +
\int^{t}_{t_0} e^{ A(t-s) } B \,dW_s
\]
and
\[
\Phi( \mathbf{w}_0 )(t) = ( e^{ A \Delta t } - I ) U_{ t_0 } +
\int^t_{t_0} e^{ A(t-s) } F( U_s ) \,ds + \int^{t}_{t_0} e^{ A(t-s) } B
\,dW_s
\]
for every $ t \in[t_0,T]$ due to (\ref{extr1a}). Since $
\operatorname{ord}( \mathbf{w}_0 ) = 1 $ (see Section
\ref{secorder}), we finally obtain
\[
U_t = e^{ A \Delta t } U_{ t_0 } + \int^t_{t_0} e^{ A(t-s) } B \,dW_s + O
( \Delta t )
\]
for the Taylor expansion corresponding to the S-wood $\mathbf{w}_0$.

%s5.1.2 ###
\subsubsection{Taylor expansion
of order $ 1 + \min(\gamma,\delta)$} \label{sec512} In order to derive
a higher order Taylor expansion, we expand the stochastic wood
$\mathbf{w}_0$. To this end, we consider the Taylor expansion given by
the S-wood $ \mathbf{w}_1 = E_{(2,1)} \mathbf{w}_0 $ (see Figure
\ref{fig2}). Here, $ \Phi( \mathbf{w}_1 ) $ and $ \Psi( \mathbf{w}_1 )
$ are presented in (\ref{extr1}) and (\ref{extr2}). Since $
\operatorname{ord}( \mathbf{w}_1 ) = 1 + \min( \gamma, \delta) $ (see
Section \ref{secorder}), we obtain
\[
U_t = e^{ A \Delta t } U_{t_0} + \biggl( \int^{\Delta t}_0 e^{ A s } \,ds \biggr) F( U_{t_0} )
+ \int^{t}_{t_0} e^{ A (t-s) } B \,dW_s + O \bigl( (\Delta t )^{ ( 1 +
\min(\gamma,\delta ) ) } \bigr)
\]
for the Taylor expansion corresponding to the S-wood $\mathbf{w}_1$.
This example corresponds to the exponential Euler scheme introduced in
\cite{jk08b}, which was already discussed in Section \ref{secsimple}
[see (\ref{onestep})].

%s5.1.3 ###
\subsubsection{Taylor expansion of order $ 1 + \min(2 \gamma, \delta) $}
\label{sec513}
For a higher order Taylor expansion, we now have several
possibilities to further expand the stochastic wood~$ \mathbf{w}_1 $.
For instance, we could consider the Taylor expansion given by the
stochastic wood $ E_{(5,1)} E_{(5,2)} \mathbf{w}_1 $ [see
(\ref{anodefig2}) for the active nodes of $\mathbf{w}_1$]. Since our
main goal is always to obtain higher order approximations with the
least possible terms and since the fifth tree of $ \mathbf{w}_1 $
[given by the nodes $ (5,1)$ and $ (5,2)$] is of order $2$ (see
Section \ref{secorder} for details), we concentrate on expanding the
lower order trees of~$ \mathbf{w}_1$. Finally, since oftentimes $
\gamma\leq\delta$ in examples (see Section \ref{spacetime} and also
Section \ref{traceclass}), we consider the stochastic wood $
\mathbf{w}_2 = E_{(4,1)} \mathbf{w}_1 $ (see Figure \ref{fig5}). It is
of order $ 1 + \min(2 \gamma, \delta) $ (see Section \ref{secorder})
and the corresponding\vspace*{1pt} Taylor approximation $ \Psi( \mathbf{w}_2 ) $ of
$ \Phi(\mathbf{w}_2) = \Delta U $ is given by $ \Psi( \mathbf{w}_2 ) =
I^0_0 + I^0_1 + I^0_2 + I^1_1[ I^0_0 ] $. This yields
\begin{eqnarray*}
U_t &=& e^{ A \Delta t } U_{t_0} + \biggl( \int^{ \Delta t }_0 e^{ A s } \,ds \biggr)
F ( U_{t_0} ) + \int^{t}_{t_0} e^{ A (t-s) }
B \,dW_s \\
&& + \int^t_{t_0} e^{ A ( t-s) } F'( U_{t_0} ) ( e^{ A \Delta s } - I
)U_{ t_0} \,ds + O \bigl( ( \Delta t )^{ ( 1 + \min(2 \gamma, \delta) ) } \bigr)
\end{eqnarray*}
for the Taylor expansion corresponding to the S-wood $\mathbf{w}_2$.

%s5.1.4 ###
\subsubsection{Taylor expansion of order $ 1 + 2 \min( \gamma, \delta) $}
\label{sec514}
In examples, we oftentimes have $ \gamma= \frac{1}{4} -
\varepsilon$, $\delta= \frac{1}{4} $ for an arbitrarily small $
\varepsilon\in(0,\frac{1}{4})$ (space--time white noise, see Section
\ref{spacetime}) or $ \gamma= \frac{1}{2} - \varepsilon$, $\delta=
\frac{1}{2} $ for an arbitrarily small $ \varepsilon\in(0,\frac{1}{2})$
(trace-class noise). In these cases, the sixth stochastic tree in $
\mathbf{w}_2 $ turns out to be the active tree of the lowest order.
Therefore, we consider the stochastic wood $ \mathbf{w}_3 = E_{(6,1)}
\mathbf{w}_2 $ (see Figure \ref{fig6}) here. It has order $ 1 + 2 \min(
\gamma, \delta) $ (see Section \ref{secorder}) and we have
\[
\Psi( \mathbf{w}_3 ) = I^0_0 + I^0_1 + I^0_2 + I^1_1[ I^0_0 ] + I^1_1[
I^0_2 ] ,
\]
which implies
\begin{eqnarray*}
U_t &=& e^{ A \Delta t } U_{t_0} + \biggl( \int^{ \Delta t }_0 e^{ A s } \,ds \biggr)
F ( U_{t_0} ) + \int^{t}_{t_0} e^{ A (t-s) } B \,dW_s
\\
&&{} +
\int^t_{t_0} e^{ A ( t-s) } F'( U_{t_0} ) ( e^{ A \Delta s } - I ) U_{
t_0} \,ds
\\
&&{} +
\int^t_{t_0} e^{ A ( t-s) } F'( U_{t_0} ) \int^s_{t_0} e^{ A ( s-r) } B
\,dW_r \,ds + O \bigl( ( \Delta t )^{ 1 + 2 \min( \gamma, \delta) } \bigr) .
\end{eqnarray*}
This example corresponds to the Taylor expansion introduced in the
beginning in Section \ref{sechigher}.

%s5.1.5 ###
\subsubsection{Taylor expansion of order
$1 + \min( 3 \gamma, \gamma+ \delta, 2 \delta)$} \label{sec515}
Since we often have $ \gamma\leq\rho$ and $ \gamma<
\frac{1}{2} $ in the examples below, the seventh stochastic tree in $
\mathbf{w}_3 $ has the lowest order in these examples. Therefore, we
consider the Taylor approximation corresponding to the S-wood $
\mathbf{w}_4 = E_{(7,1)} \mathbf{w}_3 $ (see Figure \ref{fig7}), which
is given by
\begin{eqnarray*}
U_t &=& e^{ A \Delta t } U_{t_0} + \biggl( \int^{\Delta t}_0 e^{ A s } \,ds \biggr) F ( U_{t_0}
) + \int^{t}_{t_0} e^{ A (t-s) } B \,dW_s
\\
&&{} +
\int^t_{t_0} e^{ A ( t-s) } F'( U_{t_0} ) ( e^{ A \Delta s } - I ) U_{
t_0} \,ds
\\
&&{}
+ \frac{1}{2} \int^t_{t_0} e^{ A ( t-s) } F''( U_{t_0} ) \bigl( ( e^{A
\Delta s} - I ) U_{t_0} , ( e^{A\Delta s} - I ) U_{t_0} \bigr) \,ds
\\
&&{} +
\int^t_{t_0} e^{ A ( t-s) } F'( U_{t_0} ) \int^s_{t_0} e^{ A ( s-r) } B
\,dW_r \,ds\\
&&{} + O \bigl( ( \Delta t )^{ ( 1 + \min( 3 \gamma, \gamma+ \delta, 2
\delta ) ) } \bigr) .
\end{eqnarray*}
It is of order $ 1 + \min( 3 \gamma, \gamma+ \delta, 2 \delta ) $,
which can be seen in Section \ref{secorder}.

%s5.1.6 ###
\subsubsection{Taylor expansion of order
$1 + 3 \min( \gamma, \delta, \frac{1}{3} ) $} \label{sec516}
In the case $ \gamma< \frac{1}{2} $, the 9th, 10th and 12th stochastic tree
in $\mathbf{w}_4$ all have lower or equal orders than the fifth
stochastic tree in $\mathbf{w}_4$. Therefore, we consider the S-wood $
\mathbf{w}_5 = E_{(12,1)} E_{(10,1)} E_{(9,1)} \mathbf{w}_4 $ (see
Figure \ref{fig8}) with the Taylor approximation
\begin{eqnarray*}
U_t &=& e^{ A \Delta t } U_{t_0} + \biggl( \int^{\Delta t}_0 e^{ A s } \,ds \biggr) F
( U_{t_0} ) + \int^{t}_{t_0} e^{ A (t-s) } B \,dW_s
\\&&{} +
\int^t_{t_0} e^{ A ( t-s) } F'( U_{t_0} ) ( e^{ A \Delta s } - I ) U_{
t_0} \,ds
\\&&{} +
\frac{1}{2} \int^t_{t_0} e^{ A ( t-s) } F''( U_{t_0} ) \bigl( ( e^{ A \Delta
s } - I ) U_{ t_0 }, ( e^{ A \Delta s } - I ) U_{ t_0 } \bigr) \,ds
\\&&{} +
\int^t_{t_0} e^{ A ( t-s) } F''( U_{t_0} ) \biggl( ( e^{ A \Delta s } - I )
U_{ t_0 }, \int^s_{t_0} e^{ A ( s-r) } B \,dW_r \biggr) \,ds
\\&&{} +
\frac{1}{2} \int^t_{t_0} e^{ A ( t-s) } F''( U_{t_0} ) \biggl( \int^s_{t_0}
e^{ A ( s-r) } B \,dW_r, \int^s_{t_0} e^{ A ( s-r) } B \,dW_r \biggr) \,ds
\\&&{} +
\int^t_{t_0} e^{ A ( t-s) } F'( U_{t_0} ) \int^s_{t_0} e^{ A ( s-r) } B
\,dW_r \,ds
 + O \bigl( ( \Delta t )^{ ( 1 + 3 \min( \gamma, \delta, {1}/{3}
) ) } \bigr) .
\end{eqnarray*}
\begin{remark} \label{remark1}
Not all Taylor expansions for general finite-dimensional SODEs in \cite
{kp} are used in practice due to cost and difficulty of computing the
higher iterated integrals in the expansions. For SODEs with additive
noise, however, the Wagner--Platen scheme is often used since the
iterated integrals appearing in it are linear functionals of the
Brownian motion process, thus Gaussian distributed and hence easy to
simulate. A similar situation holds for the above Taylor expansions of
SPDEs. In particular, the conditional distribution [with respect to $
F'(U_{t_0}) $] of the expression
\[
\int^t_{t_0} e^{ A ( t-s) } F'( U_{t_0} ) \int^s_{t_0} e^{ A ( s-r) } B
\,dW_r \,ds
\]
for $ t \in[t_0,T] $ in Section \ref{sec514} is Gaussian distributed
and, in principle, easy to simulate (see also Section
\ref{secnumerics}).
\end{remark}

%s5.2 ###
\subsection{Taylor expansions for
finite-dimensional SODEs}\label{SODEs}
Of course, the abstract setting for stochastic partial differential
equations of evolutionary type in Section \ref{sec2} in particular
covers the case of finite-dimensional SODEs with additive noise.
The main purpose of the Taylor expansions in this article is to
overcome the need of an It\^{o} formula in the infinite-dimensional
setting. In contrast, in the finite-dimensional case, It\^{o}'s formula
is available and the whole machinery developed here is not needed.
Nevertheless, we apply in this subsection the Taylor expansions
introduced above to stochastic ordinary differential equations with
additive noise to compare them with the well-known stochastic Taylor
expansions for SODEs in the monograph \cite{kp}.
These considerations are not so relevant in the view of applications,
since the finite-dimensional case is well studied in the literature
(see, e.g., \cite{milstein} and the above named monograph), but more
for a theoretical understanding of the new Taylor\vspace*{1pt} expansions introduced
here. More precisely, only in this subsection let $ H = \mathbb{R}^d $
with $ d \in\mathbb{N} $ be the $d$-dimensional $\mathbb{R}$-Hilbert
space of real $d$-tuples with the scalar product
\[
\langle v, w \rangle = v_1 \cdot w_1 + \cdots+ v_d \cdot w_d
\]
for every $ v = ( v_1, \ldots, v_d ) \in H $ and every $ w = ( w_1,
\ldots, w_d ) \in H $. Let also $ U = \mathbb{R}^m$ with $ m
\in\mathbb{N}$ and suppose that $ ( W_t )_{ t \in[0,T] } $ is a
standard $m$-dimensional Brownian motion. Furthermore, we suppose that
the eigenfunctions and the eigenvalues of the linear operator $ -A$ in
Assumption \ref{A1} are given by $ e_1 = ( 1, 0, \ldots, 0 ) \in H $, $
e_2 = ( 0, 1, 0, \ldots, 0 ) \in H , \ldots, e_m = ( 0, \ldots, 0, 1 )
\in H $ and $ \lambda_1 = \lambda_2 = \cdots = \lambda_d = 0 $ with $
\mathcal{I} = \{ 1, \ldots, d \}$. So, in this case $A$ is of course a
boring bounded linear operator with $ D(A) = H = \mathbb{R}^d$ and $ A
v = 0 $ for every $ v \in D(A) $. Furthermore, note that $ D ( (\kappa-
A )^r ) = H = \mathbb{R}^d $ for every $ r \in\mathbb{R} $ with an
arbitrary $ \kappa> 0 $. The bounded linear mapping $B \dvtx
\mathbb{R}^m \rightarrow\mathbb{R}^d $ is then a $d \times m$-matrix.
Due to Assumption \ref{A2}, the drift term $ F \dvtx \mathbb{R}^d
\rightarrow\mathbb{R}^d $ is then a smooth function with globally
bounded derivatives as it is assumed in \cite{kp}. The initial value $
x_0 \dvtx \Omega \rightarrow\mathbb{R}^d $ is then simply a
$\mathcal{F}_0/\mathcal{B}(\mathbb{R}^d) $-measurable mapping, which
satisfies $ \mathbb{E} | x_0 |^p < \infty $ for every $ p \in[1,\infty)
$. So, the SPDE (\ref{eqspde}) is in that case in fact a SODE and is
given by
%
%
%e23 ###
\begin{equation}\label{sode}
d U_t = F(U_t) \,dt + B \,dW_t,\qquad U_0 = u_0,
\end{equation}
for $ t \in[0,T]$.
% This SODE can also be written
% as
% $$
% U_t
% =
% u_0 +
% \int^t_0
% F( U_s )
% ds
% +
% \int^t_0 B \,dW_s
% $$
% for every $ t \in[0,T]$.
Now, we apply the abstract Taylor expansions introduced above to that
simple example. Therefore, note that the parameters in Assumption
\ref{A3} are given by $ \gamma= 1-\varepsilon$ and $ \delta=
\frac{1}{2} $ for every arbitrarily small $ \varepsilon\in(0,1)$. First
of all, we have
\[
U_t = U_{t_0} + B \cdot( W_t - W_{ t_0 } ) + O ( \Delta t )
\]
(see Section \ref{sec511}). Thus, this Taylor approximation
corresponds in the case of finite-dimensional SODEs to the Taylor
approximation for SODEs with the multi-index set
\[
\mathcal{A} = \{ v, (1), (2), \ldots, (m) \}
\]
in Theorem 5.5.1 in \cite{kp}. Here, we only mention the multi-index
set, which uniquely determines the stochastic Taylor expansion in
\cite{kp} and refer to the above named monograph for a detailed
description of the stochastic Taylor expansions for SODEs there.

The exponential Euler approximation in Section \ref{sec512} yields
%
%
%e24 ###
\begin{equation}\label{sode1}
U_t = U_{ t_0 } + F ( U_{t_0} ) \cdot( t - t_0 ) + B \cdot( W_t - W_{
t_0 } ) + O ( ( \Delta t )^{ {3/2} } ) .
\end{equation}
This is nothing else than the corresponding one-step approximation of
the classical Euler--Maruyama scheme (see Section 10.2 in \cite{kp})
and is in the setting of \cite{kp} given by the multi-index set
\[
\mathcal{A} = \{
v, (0), (1), (2), \ldots, (m)
\}
\]
in Theorem 5.5.1 there. In that sense, the name of the exponential
Euler scheme is indeed justified. While in this article, the Taylor
approximation (\ref{sode1}) is obtained via an expansion of the
$I^i_j$-operators\vspace*{1pt} (see Lemma \ref{lemformel} and Section
\ref{secintegralop}), in \cite{kp} the stochastic Taylor approximation
(\ref{sode1}) is achieved by applying It\^{o}'s formula to the
integrand $ F(U_t) $ in the SODE (\ref{sode}). Finally, the Taylor
approximation in Section~\ref{sec514} reduces to
\begin{eqnarray*}
U_t &=& U_{ t_0 } + F ( U_{t_0} ) \cdot( t - t_0 ) + B \cdot( W_t - W_{
t_0 } )
\\
&&{} + F' ( U_{ t_0 } ) \cdot B \cdot \biggl( \int^t_{t_0} \int^s_{ t_0 } dW_r \,ds
\biggr) + O ( ( \Delta t )^{ 2 } ) .
\end{eqnarray*}
The approximation above is nothing else than the one-step approximation
of the stochastic Taylor approximation given by the multi-index set
\[
\mathcal{A} = \left\{
\matrix{
v, (0), (1), (2), \ldots, (m), \cr
(1,0), (2,0), \ldots, (m,0), \cr
(1,1), (2,1), \ldots, (m,1), \cr
\vdots\cr
(1,m), (2,1), \ldots, (m,m) }
\right\}
\]
in Theorem 5.5.1 in \cite{kp}. In \cite{kp}, it is obtained via again
applying It\^{o}'s formula.

To sum up, although the method for deriving Taylor expansions in this
article ($I^i_j$-operators)\vspace*{1pt} is different to the method in
\cite{kp} (It\^{o}'s formula), the resulting Taylor approximations
above coincide.
% So, the abstract Taylor
% approximations in Hilbert spaces
% in this article
% generalize the classical
% stochastic Taylor
% expansions
% in the monograph
% \cite{kp}
% in the case of SODEs
% with additive noise
% in the sense above.

%s5.3 ###
\subsection{Simultaneous diagonalizable case}
We illustrate Assumption \ref{A3} with the case where $ A $ and $ B $
are simultaneous diagonalizable (see, for example, Section~5.5.1 in
\cite{dpz}). This assumption is commonly considered in the literature
for approximations of SPDEs (see, e.g., Section 2 in \cite{mgrb}
or see also \cite{mgrw,mgr,jk08b}). Suppose that $ U = H $ and that $B
\dvtx H \rightarrow H $ is given by
\[
B v = \sum_{i \in\mathcal{I}} b_i \langle e_i, v \rangle e_i\qquad \forall v \in H,
\]
where $ b_i $, $ i \in\mathcal{I} $, is a bounded family of real
numbers and $ e_i $, $ i \in\mathcal{I} $, is the family of
eigenfunctions of the operator $ A $ (see Assumption \ref{A1}).
Concerning Assumption \ref{A3}, note that
\begin{eqnarray*}
&&{\int^T_0} | ( \kappa- A )^{\gamma} e^{A s} B |^2_{\mathrm{HS}} \,ds
\\
&&\qquad = \sum_{i
\in\mathcal{I}} ( \kappa+ \lambda_i )^{2 \gamma} b_i^2 \biggl( \int^T_0 e^{ -
2 \lambda_i s } \,ds \biggr)
\\
&&\qquad \leq
\sum_{i \in\mathcal{I}} ( \kappa+ \lambda_i )^{2 \gamma} b_i^2 \biggl(
\int^T_0 e^{ - 2 \lambda_i s } e^{ - 2 \kappa s} \,ds \biggr) e^{ 2 \kappa T }
\\
&&\qquad =
\sum_{i \in\mathcal{I}} ( \kappa+ \lambda_i )^{2 \gamma} b_i^2\biggl (
\int^T_0 e^{ - 2 ( \kappa+ \lambda_i ) s } \,ds \biggr) e^{ 2 \kappa T }
\\
&&\qquad =
\frac{ e^{ 2 \kappa T } }{2} \biggl( \sum_{i \in\mathcal{I}} b_i^2 ( \kappa+
\lambda_i )^{ ( 2 \gamma- 1 ) } \bigl( 1 - e^{ - 2 ( \lambda_i + \kappa) T }
\bigr) \biggr)
\\
&&\qquad \leq
\frac{ e^{ 2 \kappa T } }{2} \biggl( \sum_{i \in\mathcal{I}} b_i^2 ( \kappa+
\lambda_i )^{ ( 2 \gamma- 1 ) } \biggr)
\end{eqnarray*}
for a given $ \gamma> 0 $, so
\[
{\int^T_0} | ( \kappa- A )^{\gamma} e^{A s} B |^2_{\mathrm{HS}} < \infty
\]
follows from
%
%
%e25 ###
\begin{equation} \label{diagcase}
\sum_{i \in\mathcal{I}} b_i^2 ( \kappa+ \lambda_i )^{ ( 2 \gamma- 1 ) }
< \infty
\end{equation}
for a given $ \gamma> 0 $. In this case, we also have
%
%
%e26 ###
\begin{equation} \label{diagcase2}
{\int^t_0} | e^{A s} B |^2_{\mathrm{HS}} \,ds \leq C t^{ 2 \delta}
\end{equation}
for every $ t \in[0,1] $ with $ \delta:= \min( \gamma, \frac{1}{2} ) $
and a constant $ C > 0 $.

%s5.4 ###
\subsection{Space--time white noise} \label{spacetime}

This example will be a special case of the previous one. Let $ H = L^2
( (0,1 ), \mathbb{R} ) $ be the space of equivalence classes of square
integrable measurable functions from the interval $ (0,1) $ to $
\mathbb{R} $ with the scalar product and the norm
\[
\langle u, v \rangle = \int^1_0 u(x) v(x) \,dx,\qquad | u | = \sqrt{{\int^1_0} |u(x) |^2 \,dx},
\qquad u, v \in H.
\]
Let $ U = H $ and let $ B = I \dvtx H \rightarrow H $ be the identity
operator. In addition, assume that $\alpha\dvtx (0,1)
\rightarrow\mathbb{R}$ is a bounded measurable function and that the
operator $F \dvtx H \rightarrow H $ is given by
\[
F( v )(x) := ( F v )(x) := \alpha(x) \cdot v(x),\qquad x \in( 0, 1 ),
\]
for all $ v \in H $, which clearly satisfies Assumption \ref{A2}. Also
note that
\[
F'(v) w = F(w) \quad\mbox{and}\quad F^{(i)}(v)(w_1, \ldots, w_i) = 0
\]
for all\vspace*{1pt} $ v,w,w_1, \ldots, w_i \in H $ and all $ i \in\{ 2,3, \ldots\}
$. Furthermore, let $A = \frac{\partial^{2}}{\partial x^{2}} \dvtx\break D(A)
\subset H \rightarrow H$ be the Laplace operator with Dirichlet
boundary condition, that is,
\[
A u = \sum^{\infty}_{n=1} -\lambda_n \langle e_n, u \rangle e_n,\qquad u \in H,
\]
where
\[
\lambda_n = \pi^2 n^2,\qquad e_n(x) = \sqrt{2} \sin(n \pi x ),\qquad x \in( 0, 1),
\]
for each $ n \in\mathcal{I} := \mathbb{N} $. Of course, the $ e_n $, $
n \in\mathbb{N}$, form an orthonormal basis of $ H $ (Assumption
\ref{A1}). Additionally, we choose $\kappa= 0$.

Let $ t_0 = 0 $ and $T = 1$. In view of (\ref{diagcase}), Assumption
\ref{A3} requires $\gamma= \frac{1}{4} - \varepsilon$ for every
arbitrarily small $ \varepsilon> 0 $. However, instead of
(\ref{diagcase2}), we obtain here the stronger result $ \delta=
\frac{1}{4} $, since
\[
| e^{ A s } |_{\mathrm{HS}} \leq C \biggl( \frac{1}{s} \biggr)^{1/4}
\]
for every $ s \in(0,1] $ and a constant $ C > 0 $ (see Remark 2 in
\cite{mgr}). Finally, let $ u_0 \in H $ be an arbitrary (deterministic)
function in $ H $, which satisfies Assumption \ref{A4}. The SPDE
(\ref{eqspde}) is then given by
\[
d U_t(x) = \biggl[ \frac{\partial^{2}}{\partial x^{2}} U_t(x) + \alpha(x)
U_t(x) \biggr] \,dt + d W_t,\qquad U_t(0) = U_t(1) = 0,
\]
with $U_0(x) = u_0(x)$ for $x \in(0,1)$ and $ t \in[0,1]$. After
considering Assumptions \ref{A1}--\ref{A4} for this example, we now
present the Taylor approximations in this case. Here, $
\varepsilon\in(0,\frac{1}{4})$ is always an arbitrarily small real
number in $ (0,\frac{1}{4}) $.

%s5.4.1 ###
\subsubsection{Taylor expansion of order $ 1 $}
For an approximation of $ U_t $ of order one for small $ t > 0 $, we
obtain
\[
U_t = e^{ A t } u_0 + \int^t_0 e^{ A ( t-s) } \,dW_s + O( \Delta t
)
\]
(see Section \ref{sec511}).

%s5.4.2 ###
\subsubsection{Taylor expansion
of order $ \frac{5}{4} - \varepsilon$} Here,
we have
\[
U_t = e^{ A  t } u_0 + A^{-1} ( e^{ A  t }-I ) F u_0 +
\int^t_0 e^{ A ( t-s) } \,dW_s + O \bigl( ( \Delta t )^{ ( {5/4} -
\varepsilon )} \bigr)
\]
for an approximation of order $ \frac{5}{4} - \varepsilon$ (see
Section \ref{sec512}).

%s5.4.3 ###
\subsubsection{Taylor expansion of order $\frac{5}{4}$}
In the next step, we obtain
\begin{eqnarray*}
U_t &=& e^{ A  t } u_0 + A^{-1} ( e^{ A  t }-I ) F u_0
+ \int^t_0 e^{ A ( t-s) } \,dW_s \\
&&{} + \biggl( \int^t_0 e^{ A (t-s) } F ( e^{ A  s } - I ) \,ds \biggr) u_0 + O (
( \Delta t )^{{5/4}} )
\end{eqnarray*}
for an approximation of order $ \frac{5}{4} $ (see Section
\ref{sec513}).

%s5.4.4 ###
\subsubsection{Taylor expansion
of order $\frac{3}{2} - \varepsilon$}
Here, we have
\begin{eqnarray*}
U_t &=& e^{ A  t } u_0 + A^{-1} ( e^{ A  t }-I ) F u_0 +
\int^t_0 e^{ A (t-s) } F \biggl( \int^s_0 e^{ A (s - r) } \,dW_r \biggr) \,ds
\\
&&{} + \biggl( \int^t_0 e^{ A (t-s) } F( e^{ A  s } - I ) \,ds \biggr) u_0
% \\
% &&
+ \int^t_0 e^{ A ( t-s) } \,dW_s + O \bigl( ( \Delta t )^{ ( {3/2} -
\varepsilon )} \bigr)
\end{eqnarray*}
for an approximation of order $ \frac{3}{2} - \varepsilon$ (see Section \ref{sec514}).

%s5.4.5 ###
\subsubsection{Taylor
expansion of order $\frac{7}{4}-\varepsilon$} Since $ F $ is linear
here with $F'(v) \equiv F$ and $F''(v) \equiv0$ for all $ v \in H $,
the approximation above is even more of order $\frac{7}{4} -
\varepsilon$, that is,
\begin{eqnarray*}
U_t &=& e^{ A  t } u_0 + A^{-1} ( e^{ A  t }-I ) F u_0 +
\int^t_0 e^{ A (t-s) } F \biggl( \int^s_0 e^{ A (s - r) } \,dW_r \biggr) \,ds
\\
&&{} + \biggl( \int^t_0 e^{ A (t-s) } F ( e^{ A  s } - I ) \,ds \biggr) u_0 +
\int^t_0 e^{ A ( t-s) } \,dW_s + O \bigl( ( \Delta t )^{ ( {7/4} -
\varepsilon )} \bigr)
\end{eqnarray*}
(see Section \ref{sec516}).

%s5.4.6 ###
\subsubsection{Taylor of order $ 2 - \varepsilon$}
We also consider the Taylor expansion given by the stochastic wood
\[
E_{(24,1)} E_{(22,1)} E_{(21,1)} E_{(19,1)} E_{(18,1)} E_{(16,1)}
E_{(15,1)} E_{(13,1)} \mathbf{w}_5 ,
\]
where $ \mathbf{w}_5 $ is presented in Figure \ref{fig8}. Since $ F $
is linear here, we see that the corresponding Taylor approximation is
the same as in the both examples above, so we obtain
\begin{eqnarray*}
U_t &=& e^{ A  t } u_0 + A^{-1} ( e^{ A  t }-I ) F u_0 +
\int^t_0 e^{ A (t-s) } F \biggl( \int^s_0 e^{ A (s - r) } \,dW_r \biggr) \,ds
\\
&&{} + \biggl( \int^t_0 e^{ A (t-s) } F ( e^{ A  s } - I ) \,ds \biggr) u_0 +
\int^t_0 e^{ A ( t-s) } \,dW_s + O \bigl( (\Delta t)^{( 2 - \varepsilon)} \bigr) .
\end{eqnarray*}
By further expansions, one can show that this approximation is in fact
of order $ 2 $.

%s5.5 ###
\subsection{Trace class noise} \label{traceclass}
In this subsection, we compute the smoothness parameters $ \gamma$ and
$ \delta$ in Assumption \ref{A3} for the case of trace class noise
(see, e.g., Sections~4.1 and 5.4.1 in \cite{dpz}).
This assumption is also commonly considered in the literature for
approximations of SPDEs (see, e.g., \cite{h,mgrb}). Precisely,
we suppose that $ B \dvtx U \rightarrow H $ is a Hilbert--Schmidt
operator, that is, $ | B |_{\mathrm{HS}} < \infty$. Hence, we obtain
\[
{\int^t_0} | e^{ A s } B |_{\mathrm{HS}}^2 \,ds \leq \int^t_0 ( | e^{ A s } |^2
\cdot| B |_{\mathrm{HS}}^2 ) \,ds \leq e^{2\kappa} | B |_{\mathrm{HS}}^2 t
\]
and therefore
\[
\sqrt{ {\int^t_0} | e^{ A s } B |_{\mathrm{HS}}^2 \,ds } \leq e^\kappa| B |_{\mathrm{HS}}
\sqrt{t}
\]
for all $ t \in[0,1] $. Moreover, we have
\begin{eqnarray*}
{\int^T_0} | (\kappa- A)^r e^{ A s } B |_{\mathrm{HS}}^2 \,ds &\leq& {\int^T_0} |
(\kappa- A)^r e^{ A s } |^2 | B |_{\mathrm{HS}}^2 \,ds
\\
&=& \biggl( \int^T_0 \bigl| ( \kappa- A )^r e^{ ( A - \kappa) s } e^{ \kappa s }
\bigr|^2 \,ds \biggr) | B |_{\mathrm{HS}}^2
\\
&\leq& \biggl( \int^T_0 \bigl| ( \kappa- A )^r e^{ ( A - \kappa) s } \bigr|^2 \,ds \biggr) e^{ 2
\kappa T } | B |_{\mathrm{HS}}^2
\\
&\leq& \biggl( \int^T_0 s^{ - 2 r } \,ds \biggr) e^{ 2 \kappa T } | B |_{\mathrm{HS}}^2 <
\infty
\end{eqnarray*}
for all $ r \in[0,\frac{1}{2}) $. Hence, we obtain $ \gamma=
\frac{1}{2} - \varepsilon$ and $ \delta= \frac{1}{2} $ for every
arbitrarily small $ \varepsilon\in(0,\frac{1}{2})$ in this situation.
Now, we present the Taylor expansions from Section~\ref{secabsex}
again in this special situation.

%
%s5.5.1 ###
\subsubsection{Taylor expansion of order $1$}
Here, we have
\[
U_t = e^{ A \Delta t } U_{ t_0 } + \int^t_{t_0} e^{ A(t-s) } B \,dW_s + O
( \Delta t )
\]
for a Taylor approximation of order $ 1 $ (see Section
\ref{sec511}).

%s5.5.2 ###
\subsubsection{Taylor expansion
of order $ \frac{3}{2} - \varepsilon$} \label{tcexpeuler} For a Taylor
approximation of order $ \frac{3}{2} - \varepsilon$ (see Section
\ref{sec512}), we obtain
\[
U_t = e^{ A \Delta t } U_{t_0} + \biggl( \int^{\Delta t}_{0}e^{A s} \,ds \biggr) F(
U_{t_0} ) + \int^{t}_{t_0} e^{ A (t-s) } B \,dW_s + O \bigl( ( \Delta t )^{ (
{3/2} - \varepsilon) } \bigr) .
\]
Here and below, $ \varepsilon\in(0,\frac{1}{2})$ is an arbitrarily
small real number in $ (0,\frac{1}{2}) $.

%
%s5.5.3 ###
\subsubsection{Taylor expansion of order $ \frac{3}{2} $}
The Taylor approximation in Section \ref{sec513}
reduces to
\begin{eqnarray*}
U_t &=& e^{ A \Delta t } U_{t_0} + \biggl( \int^{\Delta t}_{0}e^{A s} \,ds \biggr) F(
U_{t_0} ) + \int^{t}_{t_0} e^{ A (t-s) }
B \,dW_s \\
&&{} + \int^t_{t_0} e^{ A ( t-s) } F'( U_{t_0} ) ( e^{ A \Delta s } - I )
U_{ t_0} \,ds + O ( (\Delta t )^{ {3/2} } ) .
\end{eqnarray*}

%s5.5.4 ###
\subsubsection{Taylor expansion of
order $ 2 - \varepsilon$} Here, we obtain
\begin{eqnarray*}
U_t &=& e^{ A \Delta t } U_{t_0} + \biggl( \int^{\Delta t}_{0}e^{A s} \,ds \biggr) F(
U_{t_0} ) + \int^{t}_{t_0} e^{ A (t-s) } B
\,dW_s \\
&&{} + \int^t_{t_0} e^{ A ( t-s) } F'( U_{t_0} ) ( e^{ A \Delta s } - I )
U_{ t_0}
\,ds \\
&&{} + \int^t_{t_0} e^{ A ( t-s) } F'( U_{t_0} ) \int^s_{t_0} e^{ A (
s-r) } B \,dW_r \,ds + O \bigl( ( \Delta t )^{ ( 2 - \varepsilon) } \bigr)
\end{eqnarray*}
for a Taylor expansion of order $ 2 - \varepsilon$. This example
corresponds to the Taylor expansion introduced in Section
\ref{sec514}.

%s5.6 ###
\subsection{A special example of trace class noise}

Let $ H = U = L^2 ( ( 0, 1 )^3, \mathbb{R} ) $ be the space of
equivalence classes of square integrable measurable functions from $
(0,1 )^3$ to $ \mathbb{R} $ and consider two distinct Hilbert bases
$e_i$, $i \in\mathcal{I} := \mathbb{N}^3$, and $f_i$, $i
\in\mathcal{I}$, in $H$ given by
\[
e_i( x_1, x_2, x_3 ) = 2^{{3/2}} \sin( i_1 \pi x_1 ) \sin( i_2
\pi x_2 ) \sin( i_3 \pi x_3 )
\]
and
\begin{eqnarray*}
&&f_i( x_1, x_2, x_3 ) \\
&&\qquad=c_{(i_1 - 1)} c_{(i_2 - 1)} c_{(i_3 - 1)} \cos\bigl( (i_1-1) \pi x_1 \bigr)
\cos\bigl( (i_2-1) \pi x_2 \bigr) \cos\bigl( (i_3-1) \pi x_3 \bigr)
\end{eqnarray*}
for every\vspace*{1pt} $ i = ( i_1, i_2, i_3 ) \in\mathcal{I} = \mathbb{N}^3 $ and
every $ x = (x_1, x_2, x_3) \in (0,1 )^3 $, where $ c_n := \sqrt{2} $
for every $ n \in\mathbb{N} $ and $ c_0 = 1 $. Then, consider the
Hilbert--Schmidt operator $B \dvtx U \rightarrow H$ given by
\[
B u = \sum_{ i \in\mathbb{N}^3 } \frac{ \langle f_i , u \rangle }{ ( i_1 \cdot i_2
\cdot i_3 ) } e_i
\]
for all $ u \in U = H $. Moreover, let $ \lambda_i $, $ i
\in\mathbb{N}^3 $, be a family of real numbers given by $ \lambda_i =
\pi^2 ( i_1^2 + i_2^2 + i_3^2 ) $ for all $ i = ( i_1, i_2, i_3 )
\in\mathbb{N}^3 $. Finally,\vspace*{-2pt} consider $A = ( \frac{ \partial^2 }{
\partial x_1^2}+ \frac{ \partial^2 }{\partial x_2^2}+\frac{ \partial^2 }{
\partial x_3^2})\dvtx D(A) \subset H \rightarrow H$ (with Dirichlet boundary
conditions) given by
\[
A v = \sum_{ i \in\mathbb{N}^3 } - \lambda_i \langle e_i , v \rangle e_i
\]
for all $ v \in D(A) $, where $ D(A) $ is given by
\[
D(A) = \biggl\{ v \in H \bigg| \sum_{ i \in\mathbb{N}^3 } ( i_1^2 + i_2^2 + i_3^2
) | \langle e_i , v \rangle |^2 \biggr\} .
\]
Then, the SPDE (\ref{eqspde}) reduces to
\[
d U_t(x) = \biggl[ \biggl( \frac{ \partial^2 }{
\partial x_1^2}
+ \frac{ \partial^2 }{
\partial x_2^2}
+
\frac{ \partial^2 }{
\partial x_3^2}\biggr)
U_t(x) + F( U_t(x) ) \biggr] \,dt + \sqrt{Q} \,dW_t(x)
\]
with $ U|_{ \partial(0,1)^3 } = 0 $ for $ x \in(0,1)^3 $ and $t\in [0,T]$.
Assumptions \ref{A1}--\ref{A4} are fulfilled with $\delta= \frac{1}{2}$
and $\gamma= \frac{1}{2} - \varepsilon$ for every arbitrarily small
$\varepsilon\in(0,\frac{1}{2}) $ (see Section \ref{traceclass}). The
Taylor approximations in that situation are presented in Section
\ref{traceclass}.

%s6 ###
\section{Numerical schemes based on
the Taylor expansions}\label{secnumerics}

In this section, some numerical schemes based on the Taylor expansions
in this article are presented. We refer to \cite{j09a,j09b,jk08b,jk09c}
for estimations of the convergence orders of these schemes and also for
numerical simulations for these schemes.

For numerical approximations of SPDEs, one has to discretize both the
time interval $ [0,T] $ and the $\mathbb{R}$-Hilbert space $H$. For the
discretization of the space $H$, we use a spectral Galerkin
approximation based on the eigenfunctions of the linear operator $ A
\dvtx D(A) \subset H \rightarrow H $. More precisely, let $ (
\mathcal{I}_N )_{N \in\mathbb{N} } $ be a sequence of increasing finite
nonempty subsets of $ \mathcal{I} $, that is, $ \varnothing \neq
\mathcal{I}_N \subset \mathcal{I}_M \subset \mathcal{I} $ for all $ N,
M \in\mathbb{N} $ with $ N \leq M $ and let $ H_N :=
\operatorname{span} \langle e_i, i \in \mathcal{I}_N \rangle $ be the
finite-dimensional span of $ | \mathcal{I}_N |$-eigenfunctions for $N
\in\mathbb{N}$. The bounded linear mappings $ P_N \dvtx H \rightarrow
H_N $ are then given by $ P_N(v) = \sum_{ i \in\mathcal{I}_N } \langle e_i,
v \rangle e_i $ for every $ v \in H $ and every $ N \in\mathbb{N} $.

%s6.1 ###
\subsection{The exponential Euler scheme}

Based on the Taylor approximation in Sections \ref{secsimple} and
\ref{sec512}, we consider the family of random variables $
Y^{N,M}_k \dvtx \Omega \rightarrow H_N $, $ k=0,1,\ldots,M $, $ N, M
\in\mathbb{N}$, given by $ Y^{N,M}_0 = P_N( u_0 ) $ and
%
%
%e27 ###
\begin{eqnarray}
\label{schemeexpeuler}
Y^{N,M}_{k+1} &=& e^{ A {T}/{M} } Y^{N,M}_k +
\biggl( \int^{ {T}/{M} }_{0} e^{ A s } \,ds \biggr) ( P_N F ) ( Y^{N,M}_k )
\nonumber\\[-8pt]\\[-8pt]
&&{} + P_N \biggl( \int^{ { ( k + 1 ) T }/{ M } }_{ { k T }/{ M } } e^{ A
( {(k+1)T}/{M} - s ) } B \,dW_s \biggr)\nonumber
\end{eqnarray}
for every $k=0,1,\ldots,M-1 $ and every $N,M \in\mathbb{N}$. This
scheme is introduced and analyzed in \cite{jk08b}. As already
mentioned, it is called the \textit{exponential Euler scheme} there.

In the setting of deterministic PDEs, that is, in the case $B=0$, this
scheme reduces to
\[
Y^{N,M}_{k+1} = e^{ A {T}/{M} } Y^{N,M}_k + \biggl( \int^{ {T}/{M}
}_{0} e^{ A s } \,ds \biggr) ( P_N F ) ( Y^{N,M}_k )
\]
for every $k=0,1,\ldots,M-1 $ and every $N,M \in\mathbb{N}$. This
scheme and similar schemes, usually referred as exponential
integrators, have for deterministic PDEs been intensively studied in
the literature (see, e.g.,
\cite{bos05,bsw07,cm02,hls98,ho05a,ho05b,hos08,k05,kt05,l67,otw06}).
Such schemes are easier to simulate than may seem on the first sight
(see \cite{bsw07}). In the stochastic setting, we refer to
Sections 3 and 4 in \cite{jk08b} for a detailed description for the
simulation of the scheme (\ref{schemeexpeuler}), in particular, for the
generation of the random variables used there.

%s6.2 ###
\subsection[The Taylor scheme indicated by Section 5.1.3]{The Taylor scheme
indicated by Section \protect\ref{sec513}}

In view of Section \ref{sec513}, we obtain the Taylor scheme $
Y^{N,M}_k \dvtx \Omega \rightarrow H_N $, $ k=0,1,\ldots,M $, $ N, M
\in\mathbb{N}$, given by $ Y^{N,M}_0 = P_N( u_0 ) $ and
\begin{eqnarray*}
Y^{N,M}_{k+1} &=& e^{ A {T}/{M} } Y^{N,M}_k + \biggl( \int^{ {T}/{M}
}_{0} e^{ A s } \,ds \biggr) ( P_N F ) ( Y^{N,M}_k )
\\
&&{} + \int^{ { (k+1) T }/{ M } }_{ { k T }/{ M } } e^{ A ( {
(k+1) T }/{ M } - s ) } ( P_N F' ) ( Y^{N,M}_k )\\
&&\hspace*{61.3pt}{}\times \bigl( ( e^{ A ( s - {
k T }/{ M } ) } - I ) Y^{N,M}_k \bigr) \,ds
\\
&&{} + P_N \biggl( \int^{ { ( k + 1 ) T }/{ M } }_{ { k T }/{ M } } e^{ A
( {(k+1)T}/{M} - s ) } B \,dW_s \biggr)
\end{eqnarray*}
for every $k=0,1,\ldots,M-1 $ and every $N,M \in\mathbb{N}$.

%s6.3 ###
\subsection[The Taylor scheme indicated by
Section 5.1.4]{The Taylor scheme indicated by Section \protect\ref{sec514}}

The Taylor approximation in Section \ref{sec514} yields the Taylor
scheme $ Y^{N,M}_k \dvtx \Omega \rightarrow H_N $, $ k=0,1,\ldots,M $,
$ N, M \in\mathbb{N}$, given by $ Y^{N,M}_0 = P_N( u_0 ) $ and
\begin{eqnarray*}
Y^{N,M}_{k+1} &=& e^{ A {T}/{M} } Y^{N,M}_k + \biggl( \int^{ {T}/{M}
}_{0} e^{ A s } \,ds \biggr) ( P_N F ) ( Y^{N,M}_k )
\\
&&{} + \int^{ { (k+1) T }/{ M } }_{ { k T }/{ M } } e^{ A ( {
(k+1) T }/{ M } - s ) } ( P_N F' ) ( Y^{N,M}_k )\\
&&\hspace*{61.3pt}{}\times \bigl( ( e^{ A ( s - {
k T }/{ M } ) } - I ) Y^{N,M}_k \bigr) \,ds
\\
&&{} + \int^{ { (k+1) T }/{ M } }_{ { k T }/{ M } } e^{ A ( {
(k+1) T }/{ M } - s ) } ( P_N F' ) ( Y^{N,M}_k )\\
&&\hspace*{61.3pt}{}\times \biggl( P_N \biggl( \int^{ s }_{
{ k T }/{ M } } e^{ A ( s - u ) } B \,dW_u \biggr) \biggr) \,ds
\\
&&{} + P_N \biggl( \int^{ { ( k + 1 ) T }/{ M } }_{ { k T }/{ M } } e^{ A
( {(k+1)T}/{M} - s ) } B \,dW_s \biggr)
\end{eqnarray*}
for every $k=0,1,\ldots,M-1 $ and every $N,M \in\mathbb{N}$.

%s6.4 ###
\subsection{A Runge--Kutta scheme for SPDEs}

In principle, we can proceed with the next Taylor approximations and
obtain numerical schemes of higher order. These schemes would however
be of limited practical use due to cost and difficulty of computing the
higher iterated integrals as well as the higher order derivatives in
the Taylor approximations. Therefore, we follow a different approach
and derive a derivative free numerical scheme with simple integrals---a
so called \textit{Runge--Kutta scheme for SPDEs}. We would like to
mention that this way is the usual procedure for numerical schemes for
differential equations: Taylor expansions and their corresponding
Taylor schemes provide the underlying theory for deriving numerical
schemes, but are rarely implemented in practice. Instead of these
Taylor schemes other numerical schemes, which are easier to compute but
still depend on the Taylor expansions such as Runge--Kutta schemes or
multi-step schemes (see, e.g., \cite{db} for details) are used.

To derive a Runge--Kutta scheme for SPDEs, we consider the Taylor
approximation in Section \ref{sec514} (see also the Taylor scheme
above) from $ \frac{ k T }{ M } $ to $ \frac{ ( k + 1 ) T }{ M } $ and
obtain
\begin{eqnarray*}
U_{ (k+1) h } &\approx& e^{ A h } U_{ k h } + \biggl( \int^{ h }_{0} e^{ A s }
\,ds \biggr) F ( U_{ k h } ) + \int^{ ( k + 1 ) h }_{ k h } e^{ A ( ( k + 1 ) h
- s ) } B \,dW_s
\\&&{}
+ \int^{ (k+1) h }_{ k h } e^{ A ( (k+1) h - s ) } F' ( U_{ k h } ) \bigl( \bigl(
e^{ A ( s - k h ) } - I \bigr) U_{ k h } \bigr) \,ds
\\&&{}
+ \int^{ (k+1) h }_{ k h } e^{ A ( (k+1) h - s ) } F' ( U_{ k h } ) \biggl(
\int^s_{ k h } e^{ A ( s-r) } B \,dW_r \biggr) \,ds
\end{eqnarray*}
and hence
\begin{eqnarray*}
U_{ (k+1) h } &\approx& e^{ A h } U_{ k h } + h e^{ A h } F ( U_{ k h }
) + \int^{ ( k + 1 ) h }_{ k h } e^{ A ( ( k + 1 ) h - s ) } B \,dW_s
\\&&{}
+ h e^{ A h }
% e^{ A (
% (k+1) h - s
% ) }
F' ( U_{ k h } ) \biggl( \frac{ 1 }{ h } \int^{ (k+1) h }_{ k h } \bigl( e^{ A ( s
- k h ) } - I \bigr) U_{ k h } \,ds \biggr)
\\&&{}
+ h e^{ A h } F' ( U_{ k h } ) \biggl( \frac{ 1 }{ h } \int^{ (k+1) h }_{ k h
} \int^s_{ k h } e^{ A ( s-r) } B \,dW_r \,ds \biggr)
\\&\approx&
e^{ A h } U_{ k h } + h e^{ A h } F ( U_{ k h } ) + \int^{ ( k + 1 ) h
}_{ k h } e^{ A ( ( k + 1 ) h - s ) } B \,dW_s
\\&&{}
+ h e^{ A h }
% e^{ A (
% (k+1) h - s
% ) }
F' ( U_{ k h } ) \biggl[ \frac{ 1 }{ h } \int^{ (k+1) h }_{ k h } \bigl( e^{ A ( s
- k h ) } - I \bigr) \int^{ k h }_0 e^{ A ( k h - r ) } B \,dW_r \,ds \biggr]
\\&&{}
+ h e^{ A h } F' ( U_{ k h } ) \biggl( \frac{ 1 }{ h } \int^{ (k+1) h }_{ k h
} \int^s_{ k h } e^{ A ( s-r) } B \,dW_r \,ds \biggr)
\end{eqnarray*}
with $ h := \frac{ T }{ M } $ for $ k=0,1,\ldots,M-1$ and $M
\in\mathbb{N}$. This yields
\[
U_{ (k+1) h } \approx e^{ A h } U_{ k h } + h e^{ A h } F ( U_{ k h } +
Z^{M}_k ) + \int^{ ( k + 1 ) h }_{ k h } e^{ A ( ( k + 1 ) h - s ) } B
\,dW_s
\]
with the random variables
\begin{eqnarray*}
Z^{M}_k &=& \frac{ 1 }{ h } \int^{ (k+1) h }_{ k h } \bigl( e^{ A ( s - k h )
} - I \bigr) \int^{ k h }_0 e^{ A ( k h - r ) } B \,dW_r \,ds
\\
&&{}+ \frac{ 1 }{ h } \int^{ (k+1) h }_{ k h } \int^s_{ k h } e^{ A ( s-r)
} B \,dW_r \,ds
\end{eqnarray*}
for $k=0,1,\ldots,M-1$ and $M \in\mathbb{N}$. The corresponding
numerical scheme $ Y^{N,M}_k \dvtx \Omega\rightarrow H_N $, $
k=0,1,\ldots,M$, $N,M \in\mathbb{N}$, is then given by $ Y^{N,M}_0 =
P_N ( u_0 ) $ and
\begin{eqnarray*}
Y^{N,M}_{k+1} &=& e^{ A {T}/{M} } Y^{N,M}_k + \frac{T}{M} e^{ A
{T}/{M} } ( P_N F ) \bigl( Y^{N,M}_k + P_N ( Z^M_k ) \bigr)
\\
&&{} + P_N \biggl( \int^{ { ( k + 1 ) T }/{ M } }_{ { k T }/{ M } } e^{ A
( {(k+1)T}/{M} - s ) } B \,dW_s \biggr)
\end{eqnarray*}
for every $ k=0,1,\ldots,M-1$ and every $N,M \in\mathbb{N}$. This
Runge--Kutta scheme for SPDEs is introduced and analyzed in
\cite{j09a}. Under non-global Lipschitz coefficients of the SPDE, it is
analyzed in \cite{j09b}. Note that the random variables occurring in
the scheme above are Gaussian distributed and therefore easy to
simulate (see also Remark \ref{remark1} and \cite{j09a,j09b} for
details). More precisely, in the case of one-dimensional stochastic
reaction diffusion equations with space--time white noise it is shown in
the articles cited above that this scheme converges with the overall
order $ \frac{1}{4}$---with respect to the number of independent
standard normal distributed random variables and the number of
arithmetical operations used to compute the scheme instead of the
overall order $ \frac{1}{6} $ of classical numerical schemes (see, for
instance, \cite{g98,g99,s99,dg}) such as the linear implicit Euler
scheme.

%s7 ###
\section{Proofs}\label{secproofs}

%s7.1 ###
\subsection[Proofs of (6) and (7)]{Proofs of (\protect\ref{formel1})
and (\protect\ref{formel2})}
\begin{lemma}\label{lemformel}
Let Assumptions \ref{A1}--\ref{A4} be fulfilled and let $ i
\in\mathbb{N} $ be given. Then, we have
\[
I^0_{1^*} = I^0_1 + I^1_{1^*}[ I^0_0 ] + I^1_{1^*}[ I^0_{1^*} ] +
I^1_{1^*}[ I^0_2 ]
\]
and
\begin{eqnarray*}
I^i_{1^*}[g_1, \ldots, g_i] &=& I^i_1[g_1, \ldots, g_i]
+ I^{(i+1)}_{1^*}[ I^0_0, g_1, \ldots, g_i ] \\
&&{}
+ I^{(i+1)}_{1^*}[ I^0_{1^*}, g_1, \ldots, g_i ]
+ I^{(i+1)}_{1^*}[ I^0_2, g_1, \ldots, g_i ]
\end{eqnarray*}
for all $ g_1, \ldots, g_i \in\mathcal{P} $.
%(Equation \eqref{formel1} and equation \eqref{formel2}).
\end{lemma}
\begin{pf}
We begin with the first equation.
Since we have
\begin{eqnarray*}
F( U_s ) &=& F( U_{ t_0 } ) + \int^1_{0} F'\bigl( U_{ t_0 } + r ( U_s -
U_{t_0} ) \bigr)( U_s - U_{ t_0 } ) \,dr
\\&=&
F( U_{ t_0 } ) + \int^1_{0} F'( U_{ t_0 } + r \Delta U_s )( \Delta U_s
) \,dr
\\&=&
F( U_{ t_0 } ) + \int^1_{0} F'( U_{ t_0 } + r \Delta U_s )( I^0_0(s) )
\,dr
\\&&{}
+ \int^1_{0} F'( U_{ t_0 } + r \Delta U_s )( I^0_{1^*}(s) ) \,dr +
\int^1_{0} F'( U_{ t_0 } + r \Delta U_s )( I^0_2(s) ) \,dr
\end{eqnarray*}
for every $ s \in[t_0,T] $ due to the fundamental theorem of calculus
and (\ref{eqmain3}), we obtain
\begin{eqnarray*}
I^0_{1^*}(t) &=& \int^t_{t_0} e^{ A(t-s) } F( U_s ) \,ds
\\
&=&\int^t_{t_0} e^{ A(t-s) } F( U_{t_0} ) \,ds + \int^t_{t_0} e^{ A(t-s) }
\biggl( \int^1_{0} F'( U_{ t_0 } + r \Delta U_s )( I^0_0(s) ) \,dr \biggr) \,ds
\\
&&{}
+ \int^t_{t_0} e^{ A(t-s) } \biggl( \int^1_{0} F'( U_{ t_0 } + r \Delta U_s
)( I^0_{1^*}(s) ) \,dr \biggr) \,ds
\\
&&{}
+ \int^t_{t_0} e^{ A(t-s) } \biggl( \int^1_{0} F'( U_{ t_0 } + r \Delta U_s
)( I^0_2(s) ) \,dr \biggr) \,ds,
\end{eqnarray*}
which implies
\[
I^0_{1^*}(t) = I^0_1(t) + I^1_{1^*}[ I^0_0 ](t) + I^1_{1^*}[ I^0_{1^*}
](t) + I^1_{1^*}[ I^0_2 ](t)
\]
for all $ t \in[t_0,T] $. Moreover, we have
\begin{eqnarray*}
&& \int^1_{0} F^{(i)}( U_{ t_0 } + r \Delta U_s ) ( g_1(s), \ldots,
g_i(s) ) \frac{ (1-r)^{ (i-1) } }{ (i-1)! } \,dr
\\
&&\qquad=
\biggl[ - F^{(i)}( U_{ t_0 } + r \Delta U_s ) ( g_1(s), \ldots, g_i(s) )
\frac{ (1-r)^{i} }{ i! } \biggr]^{r=1}_{r=0}
\\
&&\qquad\quad{}
+ \int^1_0 F^{(i+1)}( U_{ t_0 } + r \Delta U_s ) ( \Delta U_s, g_1(s),
\ldots, g_i(s) ) \frac{ (1-r)^{i} }{ i! } \,dr
\\
&&\qquad=
\frac{1}{i!} F^{(i)}( U_{ t_0 } ) ( g_1(s), \ldots, g_i(s) )
\\
&&\qquad\quad{}
+ \int^1_0 F^{(i+1)}( U_{ t_0 } + r \Delta U_s ) ( \Delta U_s, g_1(s),
\ldots, g_i(s) ) \frac{ (1-r)^{i} }{ i! } \,dr
\end{eqnarray*}
for all $ s \in[t_0,T] $ and all $ g_1, \ldots, g_i \in\mathcal{P}$ due
to integration by parts and therefore, we also obtain the second
equation.
\end{pf}

%s7.2 ###
\subsection[Proof of Theorem 1]{Proof of Theorem \protect\ref{mainthm}}
For the proof of Theorem \ref{mainthm}, we need the following lemma.
%a proof of which can be found in \cite{**}.
%
\begin{lemma} \label{lem1}
Let $ X \dvtx [t_0,T] \times \Omega\rightarrow[0,\infty) $ be a
predictable stochastic process. Then, we obtain
\[
\biggl| \int^t_{t_0} X_s \,ds \biggr|_{ L^r } \leq {\int^t_{t_0}} | X_s |_{ L^r} \,ds
\]
for every $ t \in[t_0,T] $ and every $ r \in[1,\infty) $, where both
sides could be infinite.
\end{lemma}
\begin{pf}
First, we consider the case, where $ X_t \leq C $ is bounded by a
constant $ C > 0 $ for all $ t \in[0,T] $. Here, we have
\begin{eqnarray*}
\mathbb{E} \biggl[ \biggl( \int^t_{t_0} X_s \,ds \biggr)^r \biggr] &=& \int^t_{t_0} \mathbb{E} \biggl[
\biggl( \int^t_{t_0} X_u \,du \biggr)^{(r-1)}
X_s \biggr] \,ds \\
&\leq& {\int^t_{t_0}} | X_s |_{L^r} \,ds \biggl( \mathbb{E} \biggl[ \biggl( \int^t_{t_0} X_u
\,du \biggr)^{r} \biggr] \biggr)^{ (({ r-1 })/{ r }) }
\end{eqnarray*}
for every $t\in[t_0,T]$ and every $r\in[1,\infty)$ due to H\"{o}lder's inequality.
Since
\[
\mathbb{E} \biggl[ \biggl( \int^t_{t_0} X_u \,du \biggr)^{r} \biggr] < \infty
\]
is finite for every $t\in[t_0,T]$ and every $r\in[1,\infty)$ due to the boundedness of $ X \dvtx [t_0,T] \times\Omega
\rightarrow[0,\infty) $, we obtain the assertion. In the general case,
we can approximation the stochastic process $ (X_t)_{ t \in[0,T] } $ by
bounded processes $ (X^N_t)_{ t \in[0,T] } $ for $ N \in\mathbb{N} $
given by
\[
X^N_t := \min( N, X_t )
\]
for all $ t \in[0,T] $ and all $N\in\mathbb{N}$. This shows the
assertion.
\end{pf}

We also need the Burkholder--Davis--Gundy inequality in infinite
dimensions (see Lemma 7.7 in \cite{dpz}).
\begin{lemma} \label{lem2}
Let $ X\dvtx [t_0, T] \times\Omega\rightarrow \mathrm{H S} ( U, H )$ be a
predictable stochastic process with $ \mathbb{E} {\int^T_{t_0}} | X_s
|^2_{\mathrm{HS}} < \infty$. Then, we obtain
\[
\biggl| \int^t_{t_0} X_s \,dW_s \biggr|_{ L^p } \leq p \biggl( {\int^t_{t_0}} | | X_s |_{
\mathrm{HS}
} |_{L^p}^2 \,ds \biggr)^{ {1/2} }
\]
for every $ t \in[t_0,T] $ and every $ p \in[1,\infty) $, where both
sides could be infinite.
\end{lemma}

In view of the definitions of the mappings $ \Phi$ and $ \Psi$, Theorem
\ref{mainthm} immediately follows from the next lemma. For this, the
subset $ \mathbf{ST}' \subset\mathbf{ST} $ of stochastic trees given by
\begin{eqnarray*}
&&\mathbf{ST}' := \Bigl\{ \mathbf{t}= ( \mathbf{t}', \mathbf{t}'' )
\in\mathbf{ST} \big|
\forall k \in\operatorname{nd}(\mathbf{t}) \dvtx \bigl( \bigl( \exists l
\in\operatorname{nd}(\mathbf{t}) \dvtx \mathbf{t}'(l) = k \bigr)\\
&&\hspace*{168.5pt}{}\Longrightarrow \bigl( \mathbf{t}''(k) \in \{ 1, 1^* \} \bigr) \bigr) \Bigr\}
\end{eqnarray*}
is used.
\begin{lemma}\label{importantlem}
Let $ \mathbf{t}= ( \mathbf{t}', \mathbf{t}'' ) \in\mathbf{ST}' $ be an
arbitrary stochastic tree in $\mathbf{ST}'$. Then, for each $p \geq1$,
there exists a constant $C_p > 0 $ such that
\[
( \mathbb{E} [ | \Phi( \mathbf{t})(t) |^p ] )^{ {1}/{p} } \leq C_p
\cdot( t - t_0 )^{ \operatorname{ord}( \mathbf{t}) }
\]
holds for all $ t \in[t_0,T] $, where $C_p$ is independent of $ t $ and $ t_0
$ but depends on $p$, $\mathbf{t}$, $T$ and the SPDE (\ref{eqspde}).
\end{lemma}
\begin{pf}
Due to Jensen's inequality, we can assume without loss of generality
that $p \in[2,\infty) $ holds. We will prove now the assertion by induction with
respect to the number of nodes $l(\mathbf{t}) \in\mathbb{N}$.

In the base case $ l( \mathbf{t}) = 1 $, we have $ \Phi( \mathbf{t}) =
I^0_{ \mathbf{t}''(1) } $ by definition. Hence, we obtain
\begin{eqnarray*}
| \Phi( \mathbf{t})(t) |_{L^p} &=& \bigl| I^0_{ \mathbf{t}''(1) }(t) \bigr|_{L^p}
= | I^0_0(t) |_{L^p} = | ( e^{ A \Delta t } - I ) U_{ t_0 } |_{ L^p }
\\&\leq&
| ( \kappa- A )^{-\gamma} ( e^{ A \Delta t } - I ) | \cdot | ( \kappa-
A )^{\gamma} U_{ t_0 } |_{ L^p }
\\&=&
\bigl| ( \kappa- A )^{-\gamma} \bigl( e^{ ( A - \kappa) \Delta t } - e^{ -
\kappa\Delta t } \bigr) \bigr| \cdot e^{ \kappa\Delta t } \cdot | ( \kappa- A
)^{\gamma} U_{ t_0 } |_{ L^p }
\\&\leq&
\bigl| ( \kappa- A )^{-\gamma} \bigl( e^{ ( A - \kappa) \Delta t } - e^{ -
\kappa\Delta t } \bigr) \bigr| \cdot e^{ \kappa T } \cdot \Bigl( {\sup_{0 \leq s \leq
T}} | ( \kappa- A )^{\gamma} U_{ s } |_{ L^p } \Bigr)
\\&\leq&
C_p \cdot \bigl| ( \kappa- A )^{-\gamma} \bigl( e^{ ( A - \kappa) \Delta t } -
e^{ - \kappa\Delta t } \bigr) \bigr|
\end{eqnarray*}
and therefore
\begin{eqnarray*}
| \Phi( \mathbf{t})(t) |_{L^p} &\leq& C_p \cdot \bigl( \bigl| ( \kappa- A
)^{-\gamma} \bigl( e^{ ( A - \kappa) \Delta t } - I \bigr) \bigr| + | ( \kappa- A
)^{-\gamma} ( I - e^{ - \kappa\Delta t } ) | \bigr)
\\&\leq&
C_p \cdot \bigl( \bigl| ( \kappa- A )^{-\gamma} \bigl( e^{ ( A - \kappa) \Delta t } -
I \bigr) \bigr| + | ( I - e^{ - \kappa\Delta t } ) | \bigr)
\\&\leq&
C_p \cdot( \Delta t )^{ \gamma} + C_p \cdot( \Delta t ) \leq C_p \cdot(
\Delta t )^{ \gamma}
\end{eqnarray*}
for every $ t \in[t_0,T] $ in the case $ \mathbf{t}''(1) = 0 $. Here
and below, $C_p > 0$ is a constant, which changes from line to line but
is independent of $ t $ and $ t_0 $. Moreover, by Lemma \ref{lem1},
we obtain
\begin{eqnarray*}
| \Phi( \mathbf{t})(t) |_{L^p} &=& \bigl| I^0_{ \mathbf{t}''(1) }(t) \bigr|_{L^p}
= | I^0_{1^*}(t) |_{L^p} = \biggl| \int^{t}_{t_0} e^{ A(t-s) } F( U_s ) \,ds
\biggr|_{ L^p }
\\&\leq&
\int^{t}_{t_0} \bigl( \bigl| e^{ A(t-s) } F( U_s ) \bigr|_{ L^p } \bigr) \,ds \leq C_p \cdot
\biggl( {\int^{t}_{t_0}} | F( U_s ) |_{ L^p } \,ds \biggr)
\\&\leq&
C_p \cdot\biggl( \int^{t}_{t_0} ( 1 + | U_s |_{ L^p } ) \,ds \biggr) \leq C_p \cdot(
\Delta t )
\end{eqnarray*}
for every $ t \in[t_0,T] $ in the case $ \mathbf{t}''(1) = 1^* $ and
\begin{eqnarray*}
| \Phi( \mathbf{t})(t) |_{L^p} &=& \bigl| I^0_{ \mathbf{t}''(1) }(t) \bigr|_{L^p}
= | I^0_{1}(t) |_{L^p} \\
&=& \biggl| \int^{t}_{t_0} e^{ A(t-s) } F( U_{t_0} ) \,ds
\biggr|_{ L^p }
\\&\leq&
\int^{t}_{t_0} \bigl( \bigl| e^{ A(t-s) } F( U_{t_0} ) \bigr|_{ L^p } \bigr) \,ds \leq C_p
\cdot \biggl( {\int^{t}_{t_0}} | F( U_{t_0} ) |_{ L^p } \,ds \biggr)
\\&\leq&
C_p \cdot\biggl( \int^{t}_{t_0} ( 1 + | U_{t_0} |_{ L^p } ) \,ds \biggr) \leq C_p
\cdot( \Delta t )
\end{eqnarray*}
for every $ t \in[t_0,T] $ in the case $ \mathbf{t}''(1) = 1 $.
Finally, due to Lemma \ref{lem2}, we obtain
\begin{eqnarray*}
| \Phi( \mathbf{t})(t) |_{L^p} &=& \bigl| I^0_{ \mathbf{t}''(1) }(t) \bigr|_{L^p}
= | I^0_{2}(t) |_{L^p} = \biggl| \int^{t}_{t_0} e^{ A(t-s) } B \,dW_s \biggr|_{ L^p }
\\&\leq&
C_p \cdot\biggl( \int^{t}_{t_0} \bigl| \bigl| e^{ A(t-s) } B \bigr|_{\mathrm{HS}} \bigr|_{ L^p }^2 \,ds \biggr)^{
{1}/{2} }
\\&\leq&
C_p \cdot\biggl( {\int^{ ( \Delta t )}_{0} } | e^{ A s } B |^2_{\mathrm{HS}}  \,ds
\biggr)^{ {1}/{2} }
\\&\leq&
C_p \cdot( \Delta t )^{ \delta}
\end{eqnarray*}
for every $ t \in[t_0,T] $ in the case $ \mathbf{t}''(1) = 2 $. This
shows $ | \Phi( \mathbf{t})(t) |_{ L^p } \leq C_p \cdot ( \Delta t )^{
\operatorname{ord}( \mathbf{t}) } $ for every $ t \in[t_0,T]$ in the
base case $l( \mathbf{t}) = 1 $.

Suppose now that $l( \mathbf{t}) \in\{2,3,\ldots\}$. Since $
\mathbf{t}= ( \mathbf{t}', \mathbf{t}'' ) \in\mathbf{ST}' $, we must
have $ \mathbf{t}''(1) \in\{ 1, 1^* \} $. Let $ \mathbf{t}_1, \ldots,
\mathbf{t}_n \in\mathbf{ST}' $ with $ n \in\mathbb{N} $ be the subtrees
of $ \mathbf{t}$. Note that $ \mathbf{t}_1, \ldots, \mathbf{t}_n $ are
indeed in $ \mathbf{ST}' $. Then, by definition, we have
\[
\Phi( \mathbf{t})(t) = I^n_{ \mathbf{t}''(1) } [ \Phi( \mathbf{t}_1 ),
\ldots, \Phi( \mathbf{t}_n ) ](t)
\]
for every $ t \in[t_0,T] $. Therefore, by Lemma \ref{lem1}, we obtain
\begin{eqnarray*}
| \Phi( \mathbf{t})(t) |_{L^p} &=& \big|I^n_{ \mathbf{t}''(1) } [ \Phi(
\mathbf{t}_1 ), \ldots, \Phi( \mathbf{t}_n ) ](t)\big| = \big|I^n_{ 1^* } [ \Phi(
\mathbf{t}_1 ), \ldots, \Phi( \mathbf{t}_n ) ](t)\big|_{L^p}
\\
&=& \biggl| \int^{t}_{t_0} e^{ A(t-s) } \biggl( \int^1_0 F^{(n)}( U_{ t_0 } + r
\Delta U_s )\\
&&\hspace*{70.1pt}{}\times ( \Phi( \mathbf{t}_1 )(s), \ldots, \Phi( \mathbf{t}_n )(s)
) \frac{ ( 1 - r )^{n-1} }{ ( n - 1 )! } \,dr \biggr) \,ds \biggr|_{L^p}
\\&\leq&
C_p \cdot \int^{t}_{t_0} \biggl| \int^1_0 F^{(n)}( U_{ t_0 } + r \Delta U_s
)\\
&&\hspace*{51.4pt}{}\times
( \Phi( \mathbf{t}_1 )(s), \ldots, \Phi( \mathbf{t}_n )(s) ) \frac{ ( 1
- r )^{n-1} }{ ( n - 1 )! } \,dr \biggr|_{L^p} \,ds
\\&\leq&
C_p \cdot \int^{t}_{t_0} \int^1_0 \biggl| F^{(n)}( U_{ t_0 } + r \Delta U_s
)\\
&&\hspace*{51pt}{}\times
( \Phi( \mathbf{t}_1 )(s), \ldots, \Phi( \mathbf{t}_n )(s) ) \frac{ ( 1
- r )^{n-1} }{ ( n - 1 )! } \biggr|_{L^p} \,dr \,ds
\\&\leq&
C_p \cdot \int^{t}_{t_0} \int^1_0 \bigl| F^{(n)}( U_{ t_0 } + r \Delta U_s )
( \Phi( \mathbf{t}_1 )(s), \ldots, \Phi( \mathbf{t}_n )(s) ) \bigr|_{L^p} \,dr
\,ds
\\&\leq&
C_p \cdot \int^{t}_{t_0} \int^1_0 \bigl| | \Phi( \mathbf{t}_1 )(s) |
\cdots| \Phi( \mathbf{t}_n )(s) | \bigr|_{L^p} \,dr \,ds
\end{eqnarray*}
and hence
\begin{eqnarray*}
| \Phi( \mathbf{t})(t) |_{L^p} &\leq& C_p \cdot \int^{t}_{t_0} \bigl| | \Phi(
\mathbf{t}_1 )(s) | \cdots| \Phi( \mathbf{t}_n )(s) |
\bigr|_{L^p} \,ds
\\&\leq&
C_p \cdot \biggl( {\int^{t}_{t_0}} | \Phi( \mathbf{t}_1 )(s) |_{ L^{ p n } }
\cdots| \Phi( \mathbf{t}_n )(s) |_{ L^{ p n } } \,ds \biggr)
\\&\leq&
C_p \cdot \biggl( \int^{t}_{t_0} \bigl( ( \Delta s )^{ \operatorname{ord}(
\mathbf{t}_1 ) } \cdots( \Delta s )^{ \operatorname{ord}(
\mathbf{t}_n ) } \bigr) \,ds \biggr)
\\&\leq&
C_p \cdot ( \Delta t )^{ ( 1 + \operatorname{ord}(\mathbf{t}_1) +
\cdots+ \operatorname{ord}( \mathbf{t}_n ) ) } = C_p \cdot( \Delta t
)^{ \operatorname{ord}(\mathbf{t}) }
\end{eqnarray*}
for every $ t \in[t_0,T]$ in the case $ \mathbf{t}''(1) = 1^* $, since
$ l(\mathbf{t}_1), \ldots, l(\mathbf{t}_n) \leq l(\mathbf{t})-1 $ and
we can apply the induction hypothesis to the subtrees. A similar
calculation shows the result when $ \mathbf{t}''(1) = 1 $.
\end{pf}

%s7.3 ###
\subsection{Properties of the stochastic convolution} \label{secsc}
\begin{lemma} \label{lemsc}
Let Assumptions \ref{A1} and \ref{A3} be fulfilled. Then, there exists
an adapted stochastic process $O \dvtx \Omega\rightarrow C([0,T],H) $
with continuous sample paths, which is a modification of the stochastic
convolution $ \int^t_0 e^{A(t-s)} B \,dW_s $, $ t \in[0,T] $, that is, we
have
\[
\mathbb{P} \biggl[ \int^t_0 e^{A(t-s)} B \,dW_s = O_t \biggr] = 1
\]
for all $ t \in[0,T] $.
\end{lemma}
\begin{pf}
Let $ Z \dvtx [0,T] \times\Omega \rightarrow H $ be an arbitrary
adapted stochastic process with
\[
Z_t = \int^t_0 e^{ A(t-s) } B \,dW_s,\qquad \mathbb{P}\mbox{-a.s.},
\]
for all $ t \in[0,T] $. Due to Assumption \ref{A3}, such an adapted
stochastic process exists. Moreover, $ Z \dvtx [0,T]
\times\Omega\rightarrow H $ is centered and square integrable with
\[
\mathbb{E} | Z_t |^2 = \int^t_0 \bigl| e^{ A(t-s) } B \bigr|^2_{\mathrm{HS}} \,ds =
{\int^t_0}
| e^{ A s } B |^2_{\mathrm{HS}} \,ds
\]
for all $ t \in[0,T] $. Let $ \theta:= \min( \delta, \gamma) $, $ p
\geq2 $ and $ 0 \leq t_1 \leq t_2 \leq T $ be given. Then, we have
\begin{eqnarray*}
Z_{ t_2 } - Z_{t_1}
&=& \int^{ t_2 }_{ t_1 } e^{ A(t_2 - s) } B \,dW_s\\
&&{} + \int^{ t_1 }_0 \bigl( e^{ A
(t_2 - s) } - e^{ A(t_1 - s) } \bigr) B \,dW_s,\qquad
\mathbb{P}\mbox{-a.s.},
\end{eqnarray*}
and
\begin{eqnarray*}
| Z_{ t_2 } - Z_{t_1} |_{ L^p } & \leq& p \biggl( \int^{ t_2 }_{ t_1 } \bigl| \bigl|
e^{ A(t_2 - s) } B \bigr|_{\mathrm{HS}} \bigr|_{ L^p }^2 \,ds \biggr)^{ {1/2} }
\\ &&{}
+ p \biggl( \int^{ t_1 }_0 \bigl| \bigl| \bigl( e^{ A (t_2 - s) } - e^{ A(t_1 - s) } \bigr) B
\bigr|_{\mathrm{HS}} \bigr|_{L^p}^{ 2 } \,ds \biggr)^{ {1/2} }
\\ & = &
p \biggl( {\int^{ ( t_2 - t_1 ) }_{ 0 }} | e^{ A s } B |_{\mathrm{HS}}^2 \,ds
\biggr)^{ {1/2} } \\
&&{} + p \biggl( \int^{ t_1 }_0 \bigl| \bigl( e^{ A (t_2 - t_1) } - I \bigr) e^{ A s } B
\bigr|_{\mathrm{HS}}^2 \,ds \biggr)^{ {1/2} }
\end{eqnarray*}
due to Lemma \ref{lem2}.
Hence, due to Assumption \ref{A3}, we obtain
\begin{eqnarray*}
| Z_{ t_2 } - Z_{t_1} |_{ L^p } & \leq& C (t_2 - t_1)^{ \delta}\\
&&{} + C \bigl| (
\kappa- A )^{-\gamma} \bigl( e^{ A (t_2 - t_1) } - I \bigr) \bigr|
\biggl( {\int^{ t_1 }_0 }| ( \kappa- A )^{\gamma} e^{ A s } B |_{\mathrm{HS}}^2 \,ds
\biggr)^{ 1/2 }
\\ & \leq&
C (t_2 - t_1)^{ \delta} + C (t_2 - t_1)^\gamma \biggl( {\int^{ T }_0} | (
\kappa- A )^\gamma e^{ A s } B |_{\mathrm{HS}}^2 \,ds \biggr)^{ {1/2} }
\\ & \leq&
C (t_2 - t_1 )^{ \delta} + C (t_2 - t_1 )^\gamma \leq C (t_2 - t_1 )^{
\theta} ,
\end{eqnarray*}
where $ C > 0 $ is a constant changing from line to line. Since $
\theta> 0 $ is greater than zero and since $ p \geq2 $ was arbitrary,
there exists a version of $ ( Z_t )_{ t \in[0,T] } $ with continuous
sample paths due to Kolmogorov's theorem (see, e.g., Chapter 3
in \cite{dpz}).
\end{pf}

\section*{Acknowledgments}
We strongly thank the anonymous referee for his careful reading and his
very valuable advice.

% imsref loaded by lrinkeviciute, 2009-11-04 10:08:08

%
\printaddresses

\end{document}